\documentclass[review]{elsarticle}

\usepackage{amsmath,amsthm,amsfonts,amssymb,amscd,bigints}
\usepackage{subfigure}
\usepackage{graphicx,bm,color}
\usepackage{algorithm}
\usepackage{algpseudocode}
\usepackage{stmaryrd}
\usepackage{listings}
\usepackage{lineno}
\usepackage{mathtools}
\usepackage[version=3]{mhchem}
\usepackage[colorlinks,allcolors=blue]{hyperref}

\usepackage{lineno}

\definecolor{gray}{gray}{0.6}

\def\u{{\bm u}}

\def\g{{\bm g}}

\def\grad{\nabla}
\def\div{\nabla \cdot}

\begin{document}

\title{A Space Time Domain Decomposition Approach using Enhanced Velocity Mixed Finite Element Method}

\author[csm]{Gurpreet~Singh\corref{cor1}}
\ead{gurpreet@utexas.edu}
\author[csm]{Mary~F.~Wheeler}
\ead{mfw@ices.utexas.edu}

\cortext[cor1]{Corresponding author}
\address[csm]{Center for Subsurface Modeling, The University of Texas at Austin, Austin, TX 78712}

\begin{abstract}
A space-time domain decomposition approach is presented as a natural extension of the enhanced velocity mixed finite element (EVMFE) \cite{wheelerEV} for spatial domain decomposition. The proposed approach allows for different space-time discretizations on non-overlapping, subdomains by enforcing a mass continuity argument at the non-matching interface to preserve the local mass conservation property inherent to the mixed finite element methods. To this effect, we consider three different model formulations: (1) a linear single phase flow problem, (2) a non-linear slightly compressible flow and tracer transport, and (3) a non-linear slightly compressible, multiphase flow and transport. We also present a numerical solution algorithm for the proposed domain decomposition approach where a monolithic (fully coupled in space and time) system is constructed that does not require subdomain iterations. This space-time EVMFE method accurately resolves advection-diffusion transport features, in a heterogeneous medium, while circumventing non-linear solver convergence issues associated with large time-step sizes for non-linear problems. Numerical results are presented for the aforementioned, three, model formulations to demonstrate the applicability of this approach to a general class of problems in flow and transport in porous media.
\end{abstract}


\maketitle



\newcommand{\bs}[1]{\boldsymbol{#1}}


\section{Introduction}
Complex multiphase flow and reactive transport in porous medium is mathematically modeled using a system of non-linear partial differential equations. In fact, a multitude of physical processes in the porous medium are modeled by varying the nature of these non-linearities. For example, in multiphase flow the phase relative permeabilities are modeled as non-linear functions of saturation \cite{Thomas:2009tq}. On the other hand, polymer flow is modeled with an apparent viscosity term dependent on phase velocities (or phase pressure gradient) \cite{savins}. Each of these aforementioned two examples introduces a non-linear term, function of saturation or pressure, as a coefficient in the otherwise linear relationship between Darcy velocity and gradient of pressure. A consistent spatial and temporal discretization of these PDEs result in a non-linear algebraic system of equations that can be linearized using Newton method. In order to accurately resolve these non-linear terms it is often necessary that a small time-step size be used during numerical solve. For example, the timescale for reactive transport with kinetic reactions controls the flow and transport timescales \cite{Levenspiel}. This requires that small time-step sizes be used for the entire computational domain, to be consistent, wherein a majority of the reactions are restricted to only a subdomain. Consequently, these conventional approaches are often computationally prohibitive.

Several space-time domain decomposition approaches have been proposed in the past that allow non-matching spatial and temporal subdomain discretizations to address these problems.  \cite{Hughes1988, Hulbert1990} proposed space-time finite element methods for elastodynamics with discontinuous Galerkin (DG) discretization in time.  \cite{Hoang2017a} present two, space-time domain decomposition methods for a parabolic, diffusion problem using a mixed formulation differing in treatment of the non-matching subdomain interface condition. This was later extended to handle advection-diffusion type problems in \cite{Hoang2017b} using operator splitting to handle the advection and diffusion terms separately. Another space-time domain decomposition approach for parabolic problems was proposed by \cite{arbogast15}. Here, a comparison between a mortar mixed finite element method and finite volume method with ad hoc projections for coupling the non-overlapping, space-time, subdomain problems was also presented. A number of variational space-time methods for elastic wave propagation and diffusion equations were also proposed by \cite{Kocher2014, Bause2015}  using DG in space with discontinuous, continuous, and continuous differentiable Galerkin methods in time. 

In this work, we restrict ourselves to flow and transport problems in subsurface porous media with an intent to develop this approach for a more general class of problems in the near future. We extend the EVMFE method \cite{wheelerEV}, for spatial domain decomposition, to encompass the temporal dimension and propose an enhanced velocity space-time domain decomposition approach. The EVMFE method for spatial domain decomposition has been applied to a wide range of flow and transport problems \cite{Thomas:2009tq,Thomas:2011,Thomas:2012} in porous medium such as slightly and fully compressible multiphase flow as well as equation of state compositional flow. Each of these aforementioned problems differ only in the nature of the non-linearities introduced by the model formulation representing the physical process. These non-linearities occur as relative permeability and capillary pressure defined as functions of saturations, local equilibrium problem associated with phase behavior or reactive equilibrium, as well as non-linear source/sink terms in the case of reactive kinetics.

We begin by describing a linear, single-phase flow model formulation, in section \ref{sec:formulation}, with appropriate initial and boundary conditions. We then describe a mixed weak variational formulation with an appropriate choice of functional spaces for the pressure and velocity unknowns. We require that the velocity space be a subspace of $H(div)$ for local mass conservation at the non-matching space-time interface by enhancing the trace of the velocity space both spatially and temporally. With this enhanced space-time velocity space we then describe a mixed weak variational problem problem for the proposed space-time domain decomposition approach. This is followed by a fully discrete form, in section \ref{sec:discrete}, of the aforementioned weak variational form to derive the algebraic equations associated with the model formulation. A numerical solution algorithm is then presented in section \ref{sec:algo}. Here, we present a space-time monolithic solver that eliminates the need to iterate between the two subdomain problems with the interface condition serving as an evolving internal boundary condition for each of the subdomains. We then present a preliminary numerical convergence analysis for the proposed scheme in section \ref{sec:numcon}. This is followed by numerical results in section \ref{sec:results} where two practical problems of interest to flow and transport in subsurface porous media communities are considered. 

\section{Single phase flow formulation}\label{sec:formulation}
We consider the following single-phase flow model; for simplicity of discussion, represented by a first order parabolic partial differential equation system,

\begin{eqnarray}
\label{eqn:darcy}
\bs{u} = -K\nabla p \quad \text{in } \Omega\times J,\\
\label{eqn:cons}
\frac{\partial p}{\partial t} + \nabla \cdot \bs{u} = f \quad \text{in } \Omega \times J,\\
\label{eqn:bc}
p = g \quad \text{on } \partial \Omega \times J,\\
\label{eqn:init}
p = p_{0} \text{at } \Omega \times \{0\}.
\end{eqnarray}
Here, $J = (0,T]$ is the time domain of interest, $\Omega\times J = \cup_{m=1}^{q}\cup_{i=1}^{r}$ $\Omega_{i} \times J_{m}\subset \mathbb{R}^{d+1}$, d = 2 or 3, is a space-time multiblock domain in a finite $d+1$ dimensional space. Note that, in what follows, $d$ is used to indicate spatial dimensions with $d+1$ as the temporal dimension. A Dirichlet boundary condition is considered for the sake of simplicity of description and more general boundary conditions can also be treated. The subdomains $\Omega_{i}\times J_{m}$ are non-overlapping. Let us define $\Gamma^{m,n}_{i,j} = \left(\partial \Omega_{i} \times J_{m}\right) \cap \left(\partial \Omega_{j}\times J_{n}\right)$ as the interface between the space-time subdomains $\Omega_{i}\times J_{m}$ and $\Omega_{j}\times J_{n}$. Also, $\Gamma = \cup_{m,n=1}^{q}\cup_{i,j=1}^{r} \Gamma^{m,n}_{i,j}$, and $\Gamma^{m}_{i} = \left(\partial \Omega_{i}\times J_{m}\right)\cap \Gamma = \left(\partial \Omega_{i} \times J_{m}\right)\setminus \left(\Omega \times J\right)$ are the interior interfaces. Note that this description of space-time subdomain interfaces does not take into account interfaces normal to the time-direction. Since there are no arguments for pressure regularity in the time-direction, no special considerations are required. 

The functional spaces for the mixed weak formulation of Eqns. \eqref{eqn:darcy}-\eqref{eqn:init} are,\\
\begin{center}
$\bs{V} = H(div;\Omega \times J) = \left\{ \bs{v} \in \left(L^{2}(\Omega \times J)\right)^{d}: \nabla \cdot \bs{v} \in L^{2}(\Omega \times J)\right\}$,\\
$W = L^{2}(\Omega \times J)$.
\end{center}

The mixed weak formulation of Eqns. \eqref{eqn:darcy}-\eqref{eqn:init} is the pair $\bs{u}\in \bs{V}$, $p\in W$ such that,
\begin{eqnarray}
\label{eqn:wdarcy}
\left( K^{-1}\bs{u},v\right)_{\Omega\times J} = \left(p,\nabla \cdot \bs{v}\right)_{\Omega\times J} - \langle g,\bs{v}\cdot \bs{\nu}\rangle_{\partial\Omega \times J}, \quad \bs{v}\in\bs{V},\\
\label{eqn:wcons}
\left(\frac{\partial p}{\partial t}, w\right)_{\Omega\times J} + \left(\nabla \cdot \bs{u} , w\right)_{\Omega\times J} = \left(f,w\right)_{\Omega\times J}, \quad w\in W.
\end{eqnarray}

\begin{figure}[H]
\begin{center}
\includegraphics[width=7cm]{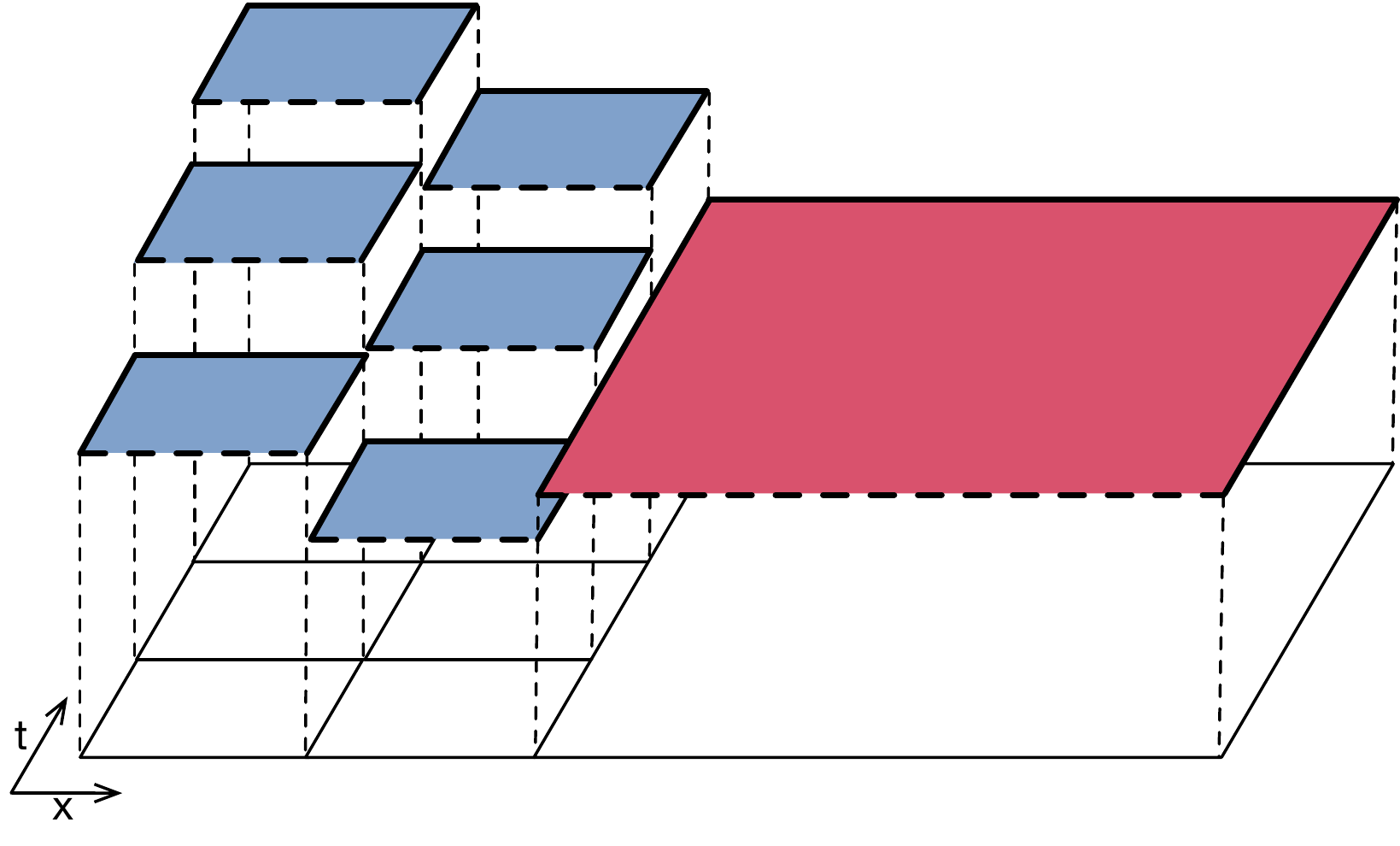}\quad
\includegraphics[width=7cm]{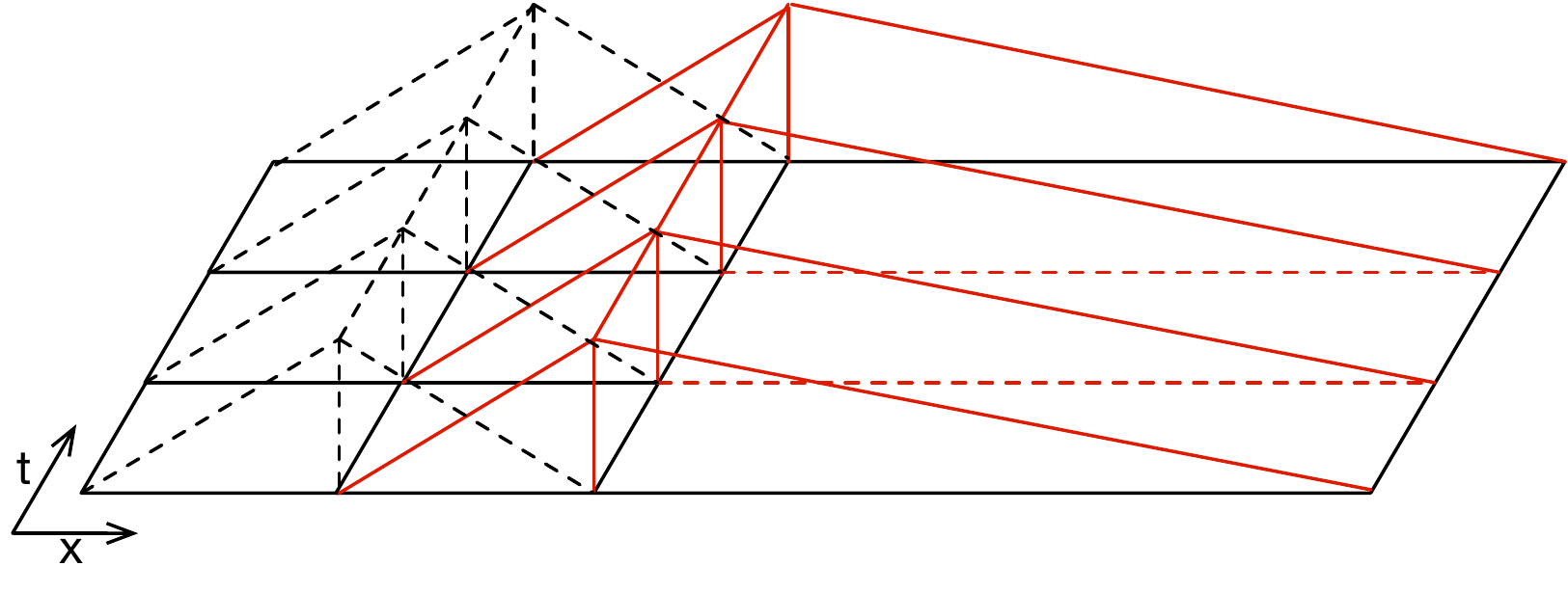}
\caption{Space-time enhanced velocity basis functions for pressure and saturation (left), and velocity (right) degrees of freedom}
\label{fig:vspace}
\end{center}
\end{figure}

Here, $\bs{\nu}$ is the normal in the spatial direction. We use the lowest order Raviart-Thomas space ($RT_{0}$) on rectangles (d = 2) or bricks (d = 3). Let $\mathcal{T}^{t,m}_{h,i}$ be a rectangular partition of $\Omega_{i}\times J_{m}$, $1\leq i \leq r$, $1\leq m \leq q$. The subdomain partitions $\mathcal{T}^{t,m}_{h,i}$ and $\mathcal{T}^{t,n}_{h,j}$ need not match on $\Gamma^{m,n}_{i,j}$. Let $\mathcal{T}^{t}_{h} =\cup_{m=1}^{q} \cup_{i=1}^{r}\mathcal{T}^{t,m}_{h,i}$. The $RT_{0}$ spaces for any element $E \in \mathcal{T}^{t}_{h}$ are defined follows:
\begin{center}
$\bs{V}^{t}_{h}(E) = \left\{ \bs{v} = (v_{1},v_{2}) \text{ or } \bs{v} = (v_{1},v_{2},v_{3}):\  v_{l} = \alpha_{l} + \beta_{l} x_{l}; \  \alpha_{l},\ \beta_{l}\  \in \  \mathbb{R}, \  l =1,...,d\right\},$\\
$W^{t}_{h}(E) = \left\{ w = \text{const}\right\}.$
\end{center} 
Note that the definition of element $E$ is in a space-time sense and hence the measure of $E$ inherits this sense accordingly. The degrees of freedom for the vector $\bs{v}\in V^{t}_{h}(E)$ can be specified by the values of the normal component $\bs{v}\cdot\bs{\nu}$ at the midpoint of all edges (faces) of E in the spatial direction. Here, $\bs{\nu}$ is the outward unit normal vector; along the spatial dimensions, on $\partial E$. The pressure finite element space on $\Omega\times J$ is defined as,\\
\begin{center}
$W^{t}_{h} = \left\{w\in L^{2}(\Omega \times J): w|_{E} \in W^{t}_{h}(E), \forall E\in \mathcal T^{t}_{h}\right\}$.\\
\end{center}
As in the original work on enhanced velocity mixed finite element method, since $W_{h}$ is a discontinuous space, no special consideration is needed on the interfaces. We want to construct a velocity finite element space $V_{h}^{t,*}\subset V$ on the multiblock partition $\mathcal{T}^{t}_{h}$ of $(\Omega\times J)$. Let, for $m=1,...,q$ and $i = 1,...,r$, 
\begin{center}
$\bs{V}^{t,m}_{h,i} = \left\{\bs{v}\in H(div;\Omega_{i}\times J_{m}):\bs{v}|_{E}\in\bs{V}_{h}(E), \forall E\in \mathcal{T}_{h,i}^{t,m}\right\}$,
\end{center}
be the usual $RT_{0}$ velocity space on $\Omega_{i}\times J_{m}$. The product space,
\begin{center}
$\bs{V}^{t}_{h} = \bs{V}^{t,1}_{h,1} \oplus \bs{V}^{t,2}_{h,1} \oplus \cdots \oplus \bs{V}^{t,1}_{h,2} \oplus \bs{V}^{t,2}_{h,2} \cdots \oplus \bs{V}^{t,q}_{h,r}$,
\end{center}
however is not a subspace of $H(div;\Omega\times J)$ since the normal vector components do not match on $\Gamma$. We therefore need to modify the degrees of freedom on $\Gamma$. Let $\mathcal{T}^{t,m,n}_{h,i,j}$ be the rectangular partition of $\mathcal{T}^{m,n}_{i,j}$ obtained from the intersection of the traces of $\mathcal{T}^{t,m}_{h,i}$ and $\mathcal{T}^{t,n}_{h,j}$. Again, as in the original work, we assume that the subdomains do not meet at an angle for a three dimensional problem. We then force the fluxes to match on each element $e \in \mathcal{T}^{t,m,n}_{h,i,j}$. Consider any element $E\in \mathcal{T}^{t,m}_{h,i}$ such that $E\cap\Gamma^{m,n}_{i,j} \ne \emptyset$. The interface grid divides the boundary edge (face) of $E$. This division can be extended inside the element in the space and time sense as shown in Figure \ref{fig:vspace} (right). 
On each subelement $E^{m}_{i}$ we define a basis function $\bs{v}_{E^{l}_{k}}$ in the $RT_{0}$ space $\bs{V}_{h}(E^{l}_{k}$) which has a normal component in each spatial dimension  $\bs{v}_{E^{l}_{k}}\nu$ equal to one on $e^{l}_{k}$ and zero on the other edges. Let $\bs{V}^{t,\Gamma}_{h}$ be the span of all such basis functions. We then define the space-time multiblock mixed finite element space to be,

\begin{center}
$\bs{V}^{t,*}_{h} = \left({\bs{V}^{t,1}_{h,1}}^{0} \oplus {\bs{V}^{t,2}_{h,1}}^{0} \oplus \cdots \oplus {\bs{V}^{t,1}_{h,2}}^{0} \oplus {\bs{V}^{t,2}_{h,2}}^{0} \cdots \oplus {\bs{V}^{t,q}_{h,r}}^{0} \oplus V_{h}^{t,\Gamma}\right) \cap H(div;\Omega \times J)$,
\end{center} 

where ${\bs{V}^{t,m}_{h,i}}^{0}$ is the subspace of $\bs{V}^{t,m}_{h,i}$ with zero normal trace on $\Gamma^{m}_{i}$. We now call $V^{t,*}_{h}$ an enhanced space-time velocity space. The additional interface degrees of freedom allows for flux continuity on the fine space-time interface grid $\mathcal{T}_{h}^{t,\Gamma} = \cup_{1\leq m \leq n}\cup_{1\leq i \leq j} \mathcal{T}^{t,m,n}_{h,i,j}$ thereby constructing an $H(div;\Omega \times J)$ - conforming velocity approximation in space and time. 

With the space-time enhanced velocity and piecewise constant pressure spaces defined, described above, the mixed finite element discretization of Eqns. \eqref{eqn:darcy}-\eqref{eqn:init} is: find $\bs{u}^{t}_{h} \in \bs{V}_{h}^{t,*}$ and $p^{t}_{h} \in W^{t}_{h}$ such that,

\begin{eqnarray}
\label{eqn:vardarcy}
\left( K^{-1}\bs{u}^{t}_{h},\bs{v}\right)_{\Omega\times J} = \left(p_{h}^{t},\nabla \cdot \bs{v}\right)_{\Omega\times J} - \langle g,\bs{v}\cdot \bs{\nu}\rangle_{\partial\Omega \times J}, \quad \bs{v}\in\bs{V}_{h}^{t,*},\\
\label{eqn:varcons}
\left(\frac{\partial p^{t}_{h}}{\partial t}, w\right)_{\Omega\times J} + \left(\nabla \cdot \bs{u}^{t}_{h} , w\right)_{\Omega\times J} = \left(f,w\right)_{\Omega\times J}, \quad w\in W_{h}^{t}.
\end{eqnarray}

\section{Fully Discrete Form}\label{sec:discrete}
In this section, we present the fully discrete form of the mixed variation formulation discussed in section \ref{sec:formulation}. We begin by describing the $RT_{0}\times DG_{0}$ basis functions as piecewise constants for pressure and saturation (or concentration) unknowns and piecewise linear in space for velocity (or flux) unknowns as follows,

\begin{equation}
w_{i}^{m} = \begin{cases}
1 & \text{on $E_{i}^{m} = x_{i-\frac{1}{2}} \leq x \leq x_{i+\frac{1}{2}} \cap t^{m} < t \leq t^{m+1}$}\\
0 & \text{otherwise}
\end{cases}
\end{equation}

\begin{equation}
\varphi^{m}_{i+\frac{1}{2}} = \begin{cases}
\frac{x-x_{i-\frac{1}{2}}}{|E_{i}^{m}|} & \text{on $E_{i}^{m}$}\\
\frac{x_{i+\frac{3}{2}}-x}{|E_{i+1}^{m}|} & \text{on $E_{i+1}^{m}$}
\end{cases}
\end{equation}
Figure \ref{fig:vspace} shows a representation these basis functions at a non-matching, space-time interface. Please note that the pressure and saturation basis functions are piecewise constants and are shown in different planes (exaggerated for clarity) to clearly indicate the jump in the time direction (dashed edge). We can obtain backward or forward Euler schemes in time by introducing the jump at the $t^{-}$ and $t^{+}$ edge, respectively of a space-time element. However, forward Euler scheme introduces stability issues requiring a CFL criteria that restricts time-step sizes. Therefore, in this work we only use the backward Euler scheme by considering a jump at the $t^{-}$ edge of the space-time element.
By construction,
\begin{equation}
\varphi_{i+\frac{1}{2}}^{m} (e^{n}_{j+\frac{1}{2}}) =
\begin{cases}
 \dfrac{1}{\rvert e^{n}_{j+\frac{1}{2}} \rvert}, \quad & i=j \text{ and } m = n,\\[3ex]
 0, \quad & \text{ otherwise}.
\end{cases}
\end{equation}
Here, $e_{j+\frac{1}{2}}^{n}$ is an edge of a space-time element.

The quadrature rules remain the same in space, as for the original enhanced velocity mixed finite element method \cite{wheelerEV}, augmented only by a midpoint rule in the time dimension as follows,

\begin{equation}
\left(v,q\right)_{TM} = \begin{cases}
 \left(v_{1},q_{1}\right)_{T\times M} & \text{if $d=1$}\\
 \left(v_{1},q_{1}\right)_{T\times M\times M} + \left(v_{2},q_{2}\right)_{M\times T \times M} & \text{if $d=2$}\\
 \left(v_{1},q_{1}\right)_{T\times M\times M \times M} + \left(v_{2},q_{2}\right)_{M\times T \times M \times M}+ \left(v_{3},q_{3}\right)_{M\times M \times T \times M} & \text{if $d=3$}
 \end{cases}
 \label{eqn:quad}
\end{equation}
As mentioned before, $ d = 1,2,3$ only represents the spatial dimension of the problem with $d+1$ as the temporal dimension. For the ease of describing the fully discrete form below below we choose one spatial dimension $d = 1$ with $d+1$ as the temporal dimension. The solution of Eqns. \eqref{eqn:darcy}-\eqref{eqn:init} can then be written as,

\begin{equation}
p = \sum_{m=1}^{q}\sum_{i = 1}^{r} P_{i}^{m} w_{i}^{m},
\end{equation}

\begin{equation}
u = \sum_{m=1}^{q}\sum_{i = 1}^{r+1} U_{i+\frac{1}{2}}^{m} \varphi_{i+\frac{1}{2}}^{m}.
\end{equation}
We will now construct an algebraic system of equations by testing the variational forms of the discrete constitutive and conservation equations with $w_{j}^{n}$ and $\varphi_{j+\frac{1}{2}}^{n}$, respectively. Here, we evaluate most of the integral terms in the variational problem on a matching grid, bifurcating to the non-matching grid only for the integrals where the non-matching, space-time interface enters the evaluation. This allows us to easily delineate the differences arising in the algebraic system of equations due to the choice of our discretization. For the constitutive equation we have,

\begin{equation}
\left( K^{-1}\bs{u},\varphi_{j+\frac{1}{2}}^{n}\right)_{\Omega\times J} = \left(p,\nabla \cdot \varphi_{j+\frac{1}{2}}^{n}\right)_{\Omega\times J} - \langle g,\varphi_{j+\frac{1}{2}}^{n}\cdot \bs{\nu}\rangle_{\partial\Omega \times J}
\label{eqn:disvardarcy}
\end{equation}

\begin{equation}
\begin{aligned}
\left( K^{-1}\bs{u},\varphi_{j+\frac{1}{2}}^{n}\right)_{\Omega\times J} & \approx  \left( K^{-1}\sum_{m=1}^{q} \sum_{i=1}^{r+1} U_{i+\frac{1}{2}}^{m}\varphi_{i+\frac{1}{2}}^{m}, \varphi_{j+\frac{1}{2}}^{n} \right)_{TM}\\
& = \frac{1}{2 \left| e_{j+\frac{1}{2}}^{n} \right|} \left( \frac{h_{j}}{K_{j}} + \frac{h_{j+1}}{K_{j+1}} \right) U_{j+\frac{1}{2}}^{n}
\end{aligned}
\end{equation}

\begin{equation}
h_{j} = x_{j+\frac{1}{2}} - x_{j-\frac{1}{2}}
\end{equation}
It is easy to see that,
\begin{equation}
\frac{\partial}{\partial x} \varphi_{i+\frac{1}{2}}^{n} = \begin{cases}
\frac{1}{\left|E_{i}^{m}\right|} & \text{on $E_{i}^{m}$}\\
\frac{1}{\left|E_{i+1}^{m}\right|} & \text{on $E_{i+1}^{m}$}
\end{cases}.
\end{equation}
The second term in Eqn. \eqref{eqn:disvardarcy} can be expanded as,
\begin{equation}
\begin{aligned}
\left(p,\nabla \cdot \varphi_{j+\frac{1}{2}}^{n}\right)_{\Omega\times J} & = \left(\sum_{m=1}^{q}\sum_{i=1}^{r} P_{i}^{m} w_{i}^{m},\nabla \cdot \varphi_{j+\frac{1}{2}}^{n}\right)_{\Omega\times J} \\
& = \bigintsss_{E_{j}^{n}} \frac{P_{j}^{n}}{\left|E_{j}^{n}\right|} + \bigintsss_{E_{j+1}^{n}} \frac{-P_{j+1}^{n}}{\left|E_{j+1}^{n}\right|}\\
& = P_{j}^{n} - P_{j+1}^{n}.
\end{aligned}
\end{equation}
Let us now consider, for a given $j_{0}\in \mathbb{I}$, a non-matching grid with fine a domain at ${\left(j_{0}+\frac{1}{2}\right)}^{-}$ and a coarse domain at ${\left(j_{0}+\frac{1}{2}\right)}^{+}$,
\begin{equation}
\begin{aligned}
\left(p,\nabla \cdot \varphi_{j_{0}+\frac{1}{2}}^{n-\frac{1}{3}}\right)_{\Omega\times J} & = \left(\sum_{m=1}^{q}\sum_{i=1}^{r} P_{i}^{m} w_{i}^{m},\nabla \cdot \varphi_{j_{0}+\frac{1}{2}}^{n-\frac{1}{3}}\right)_{\Omega\times J} \\
& = \bigintsss_{E_{j_{0}}^{n}} \frac{P_{j_{0}}^{n-\frac{1}{3}}}{\left|E_{j_{0}}^{n}\right|} + \bigintsss_{\tilde{E}_{j_{0}+1}^{n}} \frac{-P_{j_{0}+1}^{n}}{\left|\tilde{E}_{j_{0}+1}^{n}\right|}\\
& = P_{j_{0}}^{n-\frac{1}{3}} - P_{j_{0}+1}^{n}.
\end{aligned}
\label{eqn:linflu2f}
\end{equation}
Here, $\tilde{E}_{j_{0}+\frac{1}{2}}^{n}$ is a subelement of the coarse element $\tilde{E}_{j_{0}+\frac{1}{2}}^{n}$ obtained by extending the fine edge $e^{n}_{j_{0}+\frac{1}{2}}$ into the coarse element, as shown in Figure \ref{fig:vspace} (right). Similarly, testing with $\varphi_{j_{0}+\frac{1}{2}}^{n-\frac{2}{3}}$ and $\varphi_{j_{0}+\frac{1}{2}}^{n-1}$ we get,
\begin{equation}
\begin{aligned}
\left(p,\nabla \cdot \varphi_{j_{0}+\frac{1}{2}}^{n-\frac{2}{3}}\right)_{\Omega\times J} & = P_{j_{0}}^{n-\frac{2}{3}} - P_{j_{0}+1}^{n},\\
\text{and}, \\
\left(p,\nabla \cdot \varphi_{j_{0}+\frac{1}{2}}^{n-1}\right)_{\Omega\times J} & = P_{j_{0}}^{n-1} - P_{j_{0}+1}^{n},
\end{aligned}
\label{eqn:linflu2c}
\end{equation}
respectively.
Finally, the boundary or the third term in Eqn. \eqref{eqn:disvardarcy} for the coarse and fine domains is given by,
\begin{equation}
\langle g,\varphi_{j+\frac{1}{2}}^{n}\cdot \bs{\nu}\rangle_{\partial\Omega \times J} = \bigintss_{\partial E^{n}_{j+\frac{1}{2}} \cap (\partial \Omega \times J)} \frac{g_{j+\frac{1}{2}}^{n}}{\left| e_{j+\frac{1}{2}}^{n} \right|}.
\end{equation}
For the conservation equation testing with $w_{j}^{n}$ we get,
\begin{equation}
\left(\frac{\partial p}{\partial t} , w_{j}^{n}\right)_{(\Omega \times J)} - \left(\nabla \cdot u, w_{j}^{n}\right)_{(\Omega \times J)} = \left(f, w_{j}^{n}\right)_{(\Omega \times J)}.
\label{eqn:disvarcons}
\end{equation}
The first term in Eqn. \eqref{eqn:disvarcons} can be expanded as,
\begin{equation}
\begin{aligned}
\left( \frac{\partial}{\partial t} \sum_{m=1}^{q}\sum_{i=1}^{r} P_{i}^{m} w_{i}^{m}, w_{j}^{n} \right)_{\Omega \times J}  & =\left( \frac{\partial}{\partial t} P_{j}^{n} w_{j}^{n} , w_{j}^{n}\right) + \left( P_{j}^{n}-P_{j}^{n-1}, w_{j}^{n-1}\right)\\
& = \left( P_{j}^{n} - P_{j}^{n-1}\right) |E_{j}^{n-1}|.
\end{aligned}
\end{equation}
For the coarse domain $(j_{0}+\frac{1}{2})^{+}$, the above integral remains unchanged. However, for the fine domain $(j_{0}+\frac{1}{2})^{-}$ we have,
\begin{equation}
\begin{aligned}
\left(\frac{\partial p}{\partial t} , w_{j_{0}}^{n-\frac{2}{3}}\right)_{(\Omega \times J)} = & \left( P_{j_{0}}^{n-\frac{2}{3}} - P_{j_{0}}^{n-1}\right) |E_{j_{0}}^{n-1}|, \\
\left(\frac{\partial p}{\partial t} , w_{j_{0}}^{n-\frac{1}{3}}\right)_{(\Omega \times J)} = & \left( P_{j_{0}}^{n-\frac{1}{3}} - P_{j_{0}}^{n-\frac{2}{3}}\right) |E_{j_{0}}^{n-\frac{2}{3}}|, \\
\left(\frac{\partial p}{\partial t} , w_{j_{0}}^{n}\right)_{(\Omega \times J)} = & \left( P_{j_{0}}^{n} - P_{j_{0}}^{n-\frac{1}{3}}\right) |E_{j_{0}}^{n-\frac{1}{3}}|.
\end{aligned}
\end{equation}
The second term in Eqn. \eqref{eqn:disvarcons} is, 
\begin{equation}
\begin{aligned}
\left(\nabla \cdot \bs{u}, w_{j}^{n}\right)_{\Omega \times J} & = \left(\nabla \cdot \bs{u}, w_{j}^{n}\right)_{E_{j}^{n}} \\
& = \int_{\partial E_{j}^{n}}  \bs{u} \cdot \nu_{E_{j}^{n}}\\
& = \int_{e_{j+\frac{1}{2}}^{n}} \bs{u}  - \int_{e_{j-\frac{1}{2}}^{n}} \bs{u} \\
& = \bigintsss_{e_{j+\frac{1}{2}}^{n}} \sum_{m=1}^{q}\sum_{i = 1}^{r+1} U_{i+\frac{1}{2}}^{m} \varphi_{i+\frac{1}{2}}^{m}  - \bigintsss_{e_{j-\frac{1}{2}}^{n}} \sum_{m=1}^{q}\sum_{i = 1}^{r+1} U_{i+\frac{1}{2}}^{m} \varphi_{i+\frac{1}{2}}^{m} \\
& = \bigintsss_{e_{j+\frac{1}{2}}^{n}} \frac{U_{j+\frac{1}{2}}^{n}}{\left|e_{j+\frac{1}{2}}^{n}\right|}  -\bigintsss_{e_{j-\frac{1}{2}}^{n}} \frac{U_{j-\frac{1}{2}}^{n}}{\left|e_{j-\frac{1}{2}}^{n}\right|} \\ 
& = U_{j+\frac{1}{2}}^{n} - U_{j-\frac{1}{2}}^{n}.
\end{aligned}
\end{equation}
For a fine domain element, with an edge at the the non-matching space time interface, we can write this term as,
\begin{equation}
\begin{aligned}
\left(\nabla \cdot \bs{u}, w_{j_{0}}^{n-\frac{2}{3}}\right)_{\Omega \times J} & = U_{j_{0}+\frac{1}{2}}^{n-\frac{2}{3}} - U_{j_{0}-\frac{1}{2}}^{n-\frac{2}{3}}\\
\left(\nabla \cdot \bs{u}, w_{j_{0}}^{n-\frac{1}{3}}\right)_{\Omega \times J} & = U_{j_{0}+\frac{1}{2}}^{n-\frac{1}{3}} - U_{j_{0}-\frac{1}{2}}^{n-\frac{1}{3}}\\
\left(\nabla \cdot \bs{u}, w_{j_{0}}^{n}\right)_{\Omega \times J} & = U_{j_{0}+\frac{1}{2}}^{n} - U_{j_{0}-\frac{1}{2}}^{n}
\end{aligned}
\label{eqn:lincon2f}
\end{equation}
Similarly, for a coarse domain element,
\begin{equation}
\begin{aligned}
\left(\nabla \cdot \bs{u}, w_{j_{0}+1}^{n}\right)_{\Omega \times J} & = U_{j_{0}+\frac{3}{2}}^{n} - U_{j_{0}+\frac{1}{2}}^{n-\frac{2}{3}} - U_{j_{0}+\frac{1}{2}}^{n-\frac{1}{3}} - U_{j_{0}+\frac{1}{2}}^{n}\\
\end{aligned}
\label{eqn:lincon2c}
\end{equation}
The third term in Eqn. \eqref{eqn:disvarcons} is easy to evaluate for both the coarse and fine domains as,
\begin{equation}
\begin{aligned}
\left(f, w_{j}^{n}\right) _{\Omega \times J} = f_{j}^{n} \left|E_{j}^{n}\right|.
\end{aligned}
\label{eqn:linconsource}
\end{equation}
We now have a system of algebraic equations from the fully discrete form of the variational problem \eqref{eqn:vardarcy}-\eqref{eqn:varcons}. For a non-linear, multiphase flow problems the above fully discrete form results in a non-linear, algebraic, system of equations.

\section{Solution Algorithm}\label{sec:algo}
\begin{figure}[H]
\begin{center}
\includegraphics[width=7cm]{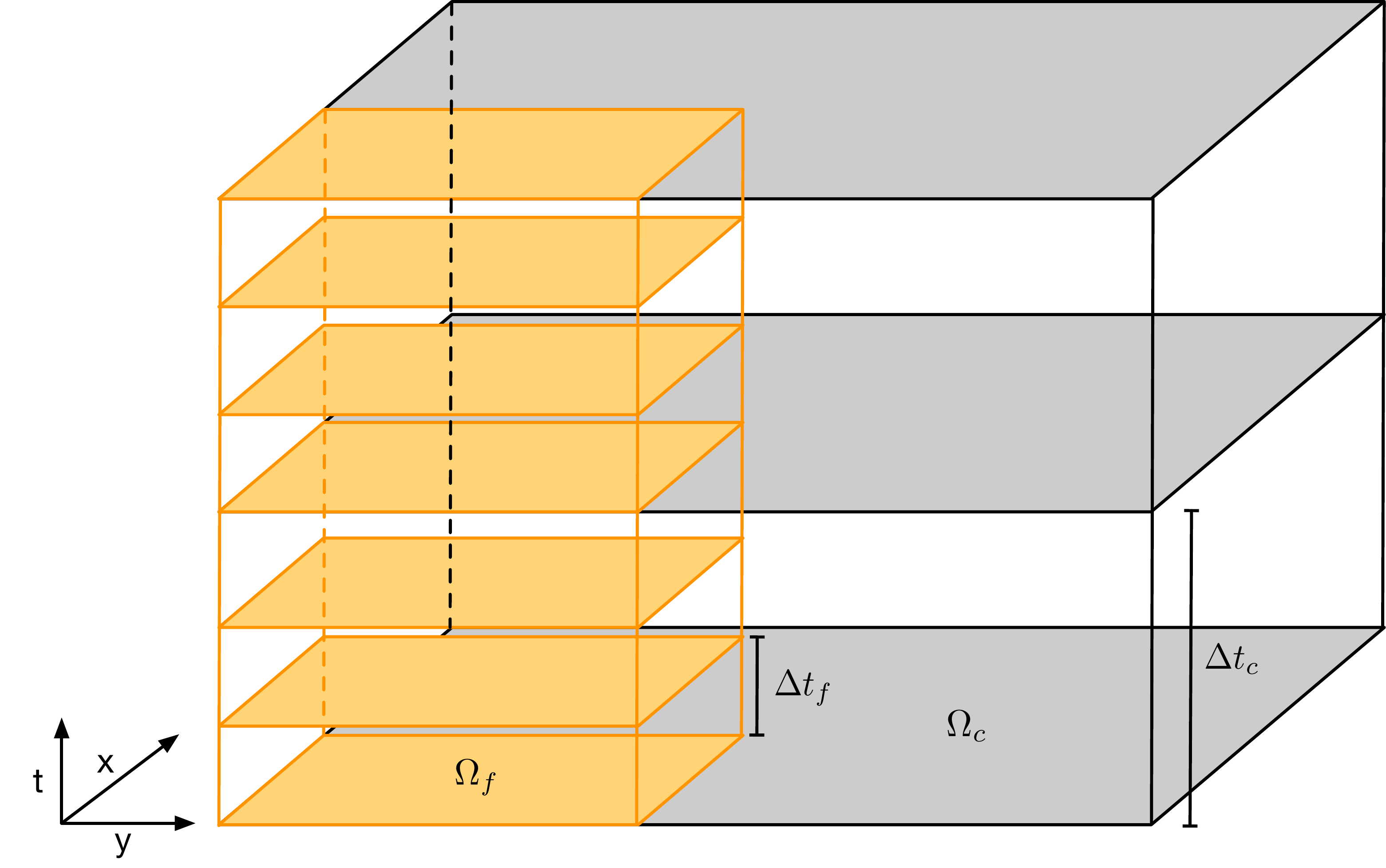}\quad
\caption{Space-time monolithic system construction schematic}
\label{fig:scheme}
\end{center}
\end{figure}

We propose a space-time monolithic solver for the algebraic system obtained in the previous section. Here we first define define matching times as the location along the time dimension where the coarse and fine domain boundaries ($\partial \Omega_{c} \cap \partial \Omega_{f}$) overlap such that $ \frac{\Delta t_{c}}{\Delta t_{f}} = l \in \mathbb I $, where $l$ is a non-zero integer. Figure \ref{fig:scheme} shows a schematic of the space-time domain decomposition used for the construction of the monolithic system with $l = 3$. 
Although the single phase flow model Eqns. \eqref{eqn:darcy}-\eqref{eqn:init} are linear, we still describe a Newton linearization step to generalize the solution algorithm for extensions to non-linear multiphase flow and reactive transport problems in sub-surface porous media. The numerical results section presents experiments where the proposed space-time domain decomposition scheme is used for two non-linear flow and transport problems. A monolithic system is then constructed to solve the algebraic equations in the space-time unknowns over the coarse and fine domains for each matching time-step ($\Delta t_{c}$) as follows,
\[
\left[ 
\begin{array}{c@{}c@{}c@{}c@{}c}
 & \vdots & \vdots & \vdots &  \\
  & \textcolor{red}{n-\frac{1}{3}} & \textcolor{red}{n-\frac{2}{3}} & \textcolor{red}{n}&  \\
   & \textcolor{red}{temporal} & spatial & &  \\
  \hdots & \left[\begin{array}{ccc}
                      \times&0&0\\ 
                      0&\times&0\\
                      0&0&\times\\
                      \end{array}\right] 
                      & \left[\begin{array}{ccc}
                      \times  &  \times & 0\\ 
                      \times  & \times & \times \\
                      0 & \times & \times \\
                      \end{array}\right]
                      &\textbf{0}
                      & \hdots \\    
& & \textcolor{red}{temporal} & spatial &   \\     
  \hdots & \textbf{0}
                      & \left[\begin{array}{ccc}
                      \times  &  0 & 0\\ 
                     0  & \times  & 0 \\
                      0 & 0 & \times \\
                      \end{array}\right]
                      & \left[\begin{array}{ccc}
                      \times &  \times & 0\\ 
                      \times & \times & \times\\
                      0 & \times & \times \\
                      \end{array}\right]
                      & \hdots \\
                       & \vdots & \vdots & \vdots &
                       \end{array}\right]
\left[
\begin{array}{c@{}}
\vdots\\
\delta p_{i-1}^{n-\frac{1}{3}}\\
\delta p_{i}^{n-\frac{1}{3}}\\
\delta p_{i+1}^{n-\frac{1}{3}}\\
\delta p_{i-1}^{n-\frac{2}{3}}\\
\delta p_{i}^{n-\frac{2}{3}}\\
\delta p_{i+1}^{n-\frac{2}{3}}\\
\delta p_{i-1}^{n}\\
\delta p_{i}^{n}\\
\delta p_{i+1}^{n}\\
\vdots
\end{array}
\right]
=
\left[
\begin{array}{c@{}}
\vdots\\
-R_{i-1}^{n-\frac{1}{3}}\\
-R_{i}^{n-\frac{1}{3}}\\
-R_{i+1}^{n-\frac{1}{3}}\\
-R_{i-1}^{n-\frac{2}{3}}\\
-R_{i}^{n-\frac{2}{3}}\\
-R_{i+1}^{n-\frac{2}{3}}\\
-R_{i-1}^{n}\\
-R_{i}^{n}\\
-R_{i+1}^{n}\\
\vdots
\end{array}
\right].
\]  
Note that here the flux unknowns $\delta u$ are eliminated by taking a Schur-complement of the original linear algebraic system in pressure and flux unknowns $\delta p$ and $\delta u$, respectively. We therefore avoid the saddle-point system associated with the original, linear algebraic system. Although not restrictive, for the ease of description we assumed non-matching grids in the time-dimension only with matching grids in space for the coarse and fine domains. The spatial sub-matrix then has a known sparsity pattern of three, five, and seven non-zero diagonals for one, two, and three spatial dimensions, respectively for RT0 mixed finite element discretization in space. As for the original EVMFEM scheme in space, the spatial sparsity pattern alters when non-matching spatial grids are considered. The temporal sub-matrix is always diagonal for DG0 discretization in time. 

\begin{figure}[h]
\begin{center}
\includegraphics[width=7cm]{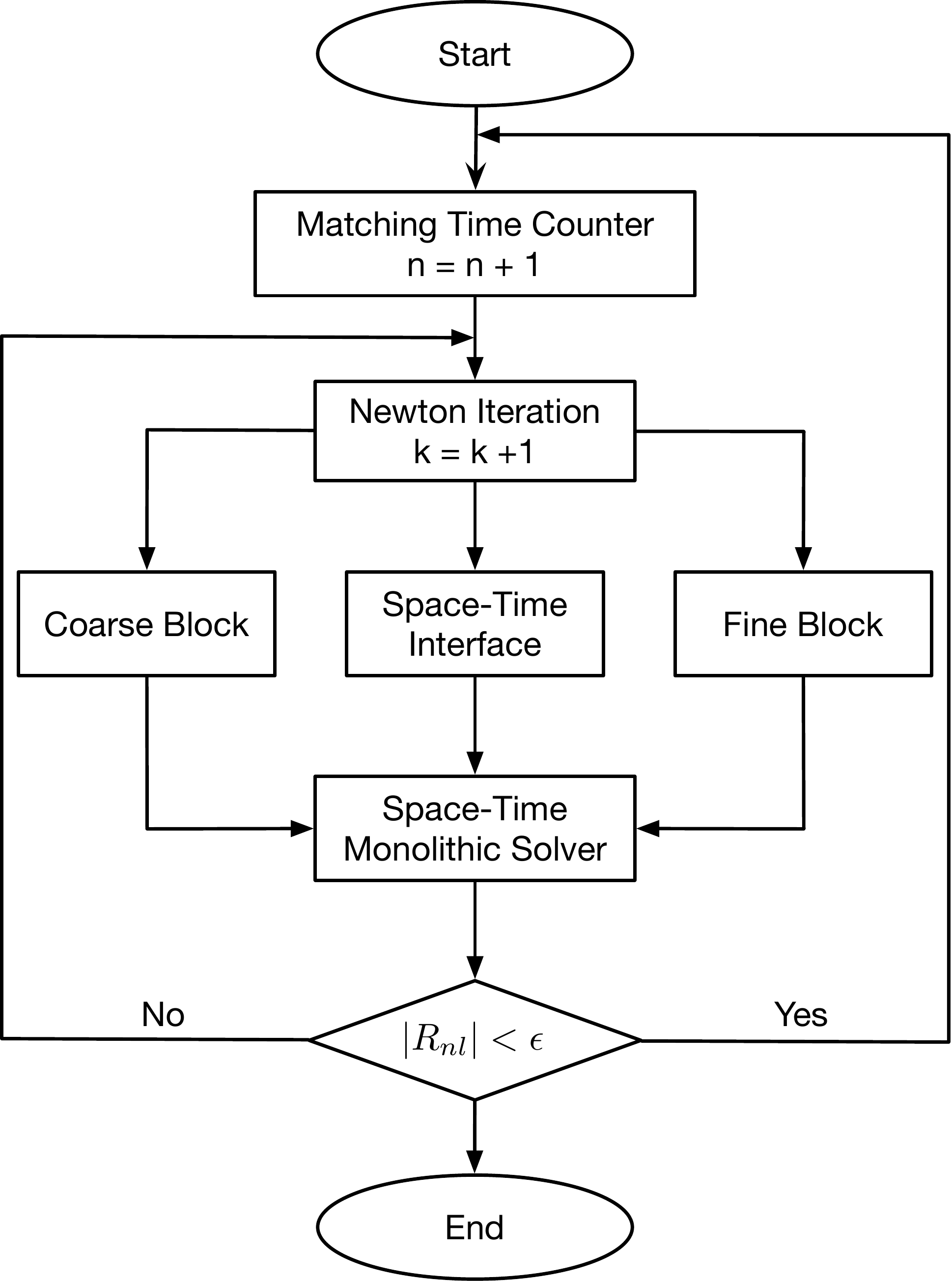}
\caption{Solution algorithm for space-time monolithic system construction and non-linear solver}
\label{fig:flow}
\end{center}
\end{figure}

Since DG0 in time is closely related to the backward Euler scheme \cite{arbogast15}, the solution scheme presented here is fully-implicit in the space-time unknowns. Figure \ref{fig:flow} shows this fully-implicit solution algorithm for the space-time domain decomposition approach proposed in this work. Here, $n$ and $k$ are the time-step and Newton iteration counters, respectively. For the single phase flow problem \eqref{eqn:darcy}-\eqref{eqn:init}, a single Newton iteration (k=1) is required for the linear system to converge, as expected. For non-linear multiphase flow and transport problems, we require that the max norm of the non-linear $\lvert R_{nl} \rvert$ residuals be less than a desired tolerance $\epsilon$. We initially rely upon direct-solvers for the purpose of testing and benchmarking the solution algorithm. However, the monolithic space-time solver allows us to utilize parallel in time linear solvers and preconditioners presented in \cite{falgout14} with relative ease. This also renders us a massively parallel, time-concurrent, framework for solving general, sub-surface, non-linear flow and reactive transport problems.

\section{Numerical convergence analysis}\label{sec:numcon}
\begin{figure}[H]
\begin{center}
\includegraphics[width=5cm,trim=10cm 2cm 9cm 1cm, clip]{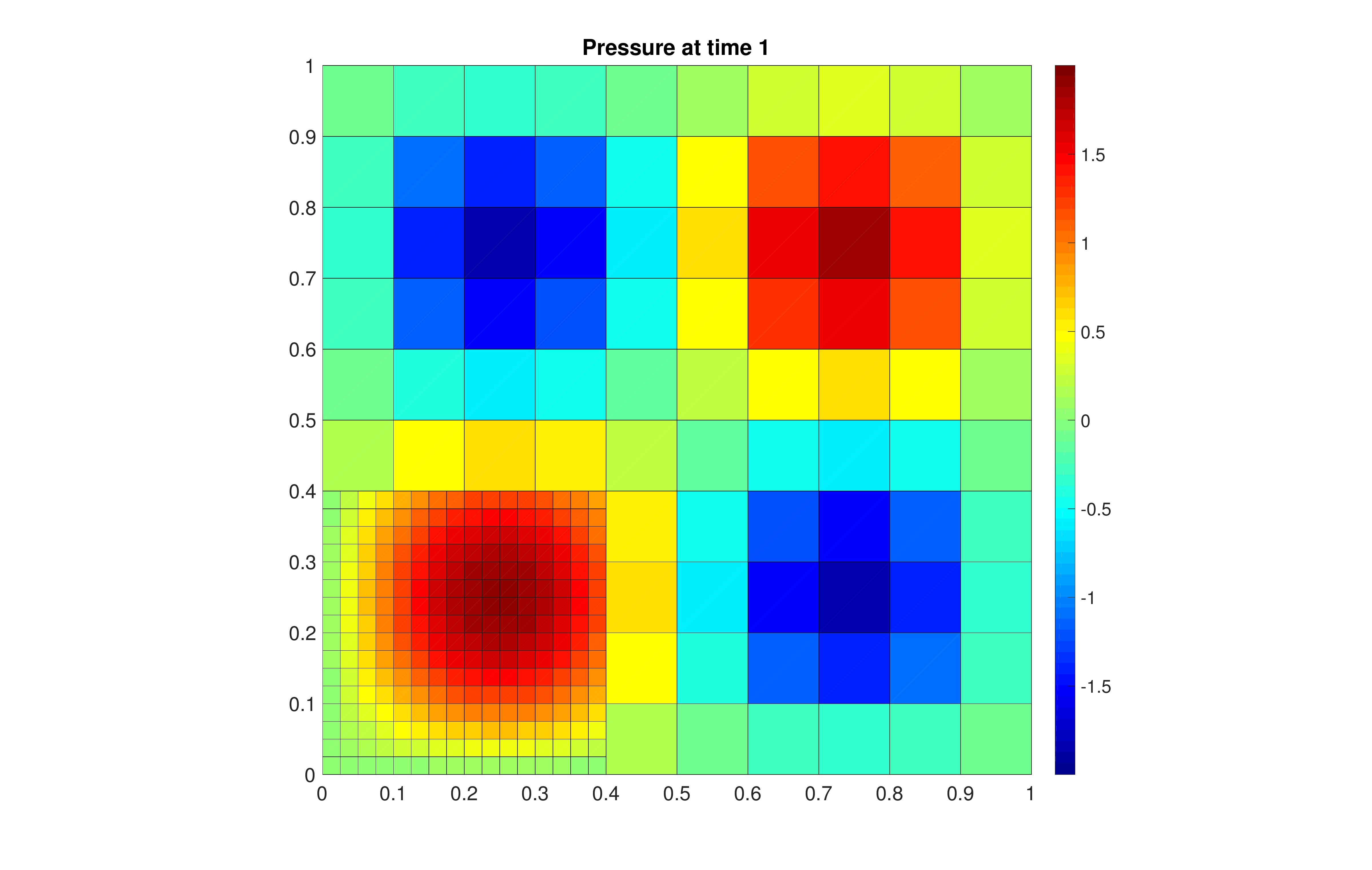}
\includegraphics[width=5cm,trim=10cm 2cm 9cm 1cm, clip]{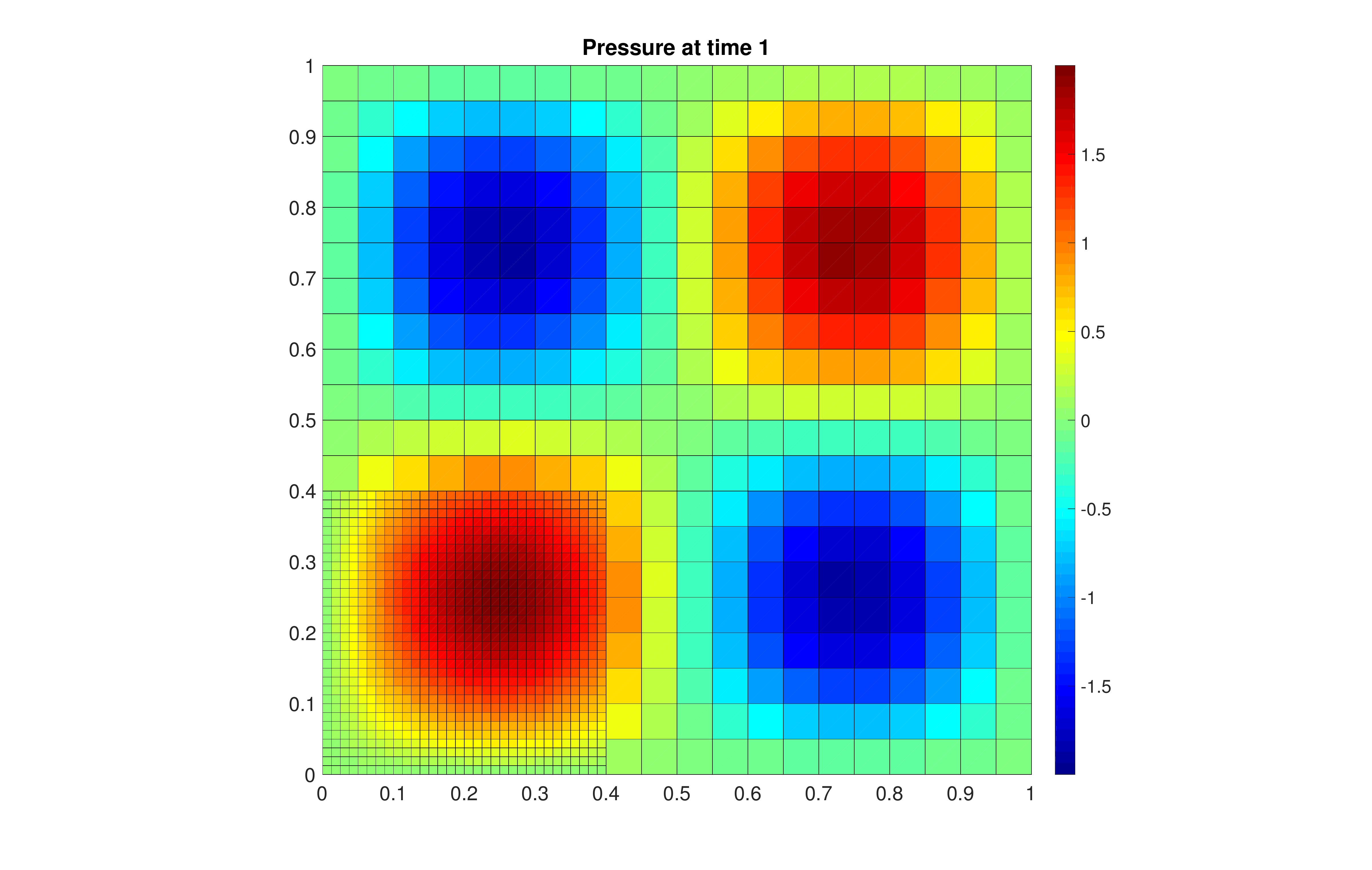}
\includegraphics[width=5cm,trim=10cm 2cm 9cm 1cm, clip]{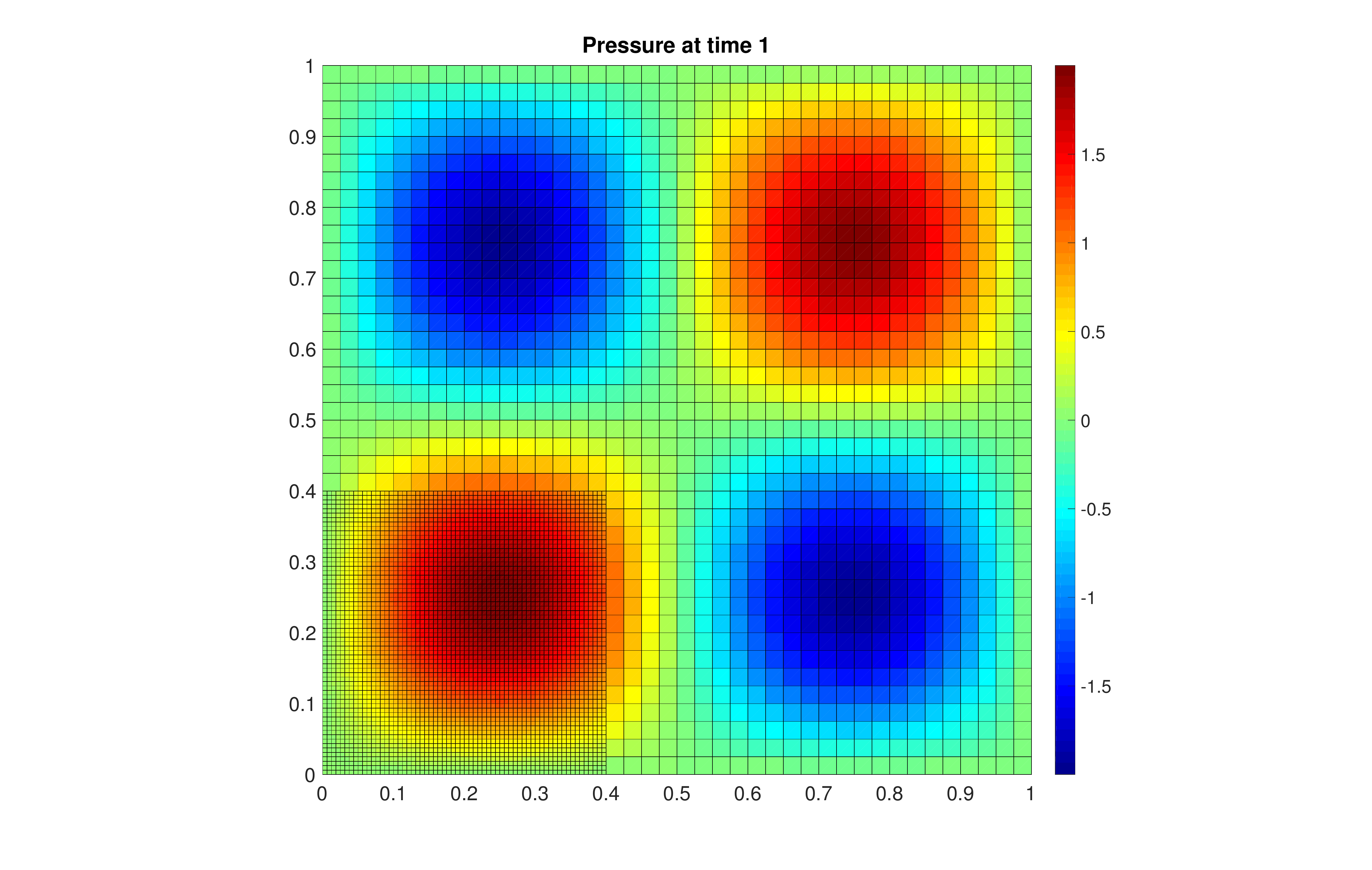}
\caption{Pressure solution for the three refinement levels}
\label{fig:manu}
\end{center}
\end{figure}

In this section, we consider a linear parabolic problem to verify numerical convergence of the proposed space-time domain decomposition approach as described below,
\begin{equation}
\frac{\partial(p)}{\partial t} + \nabla \cdot \bs{u} = f  \text{~in~} \Omega \times J,
\end{equation}
\begin{equation}
\u = -I(\nabla p) \text{~in~} \Omega \times J.
\end{equation}
Here, $I$ is the identity matrix with other symbols with their usual meanings. Further, the boundary and initial conditions are chosen to be,
\begin{equation}
p = 0 \quad\text{~on~} \partial \Omega \times J,
\end{equation}
and
\begin{equation}
p = 1 \text{~at~} \Omega \times \{0\}, 
\end{equation}
respectively. The following choice of the forcing function,
\begin{equation}
f = e^{c_{1}t}\left(c_{1} + 8\pi ^{2}\right)sin(2\pi x) sin(2 \pi y),
\end{equation}
then allows us to obtain an analytical form of the solution as,
\begin{equation}
p = e^{c_{1} t}sin(2\pi x) sin(2 \pi y).
\end{equation}
An $L^{2}$ norm in space and time for the error in pressure $||p-p_{h}||_{L^{2}}$ is defined as follows,
\begin{equation}
||p-p_{h}||_{\Omega\times J} = \left(\int_{\Omega\times J} (p-p_{h}^{\tau})^{2}\right)^{1/2} = \sum_{i}^{q}\sum_{m}^{r} \left((p-p_{i}^{m})^{2}\lvert E_{i}^{m}\rvert\right)^{1/2}
\end{equation}
We select a unit domain with two spatial (d = 2) and one temporal (d+1 = 3) dimension of size 1 $\times$ 1 $\times$ 1. The fine subdomain is refined by a factor of 4 with respect to the coarse subdomain, both spatially and temporally. Equipped with the above norm, we compute the errors in the coarse and fine subdomains for increasing mesh refinement while maintaining the aforementioned refinement factor of 4 between the coarse and fine subdomains. Here, we only show a sense of convergence for the proposed scheme and reserve optimal convergence rate studies using appropriate norms for a future work along with derivation of rigorous a-priori error estimates. 
\begin {table}[H]
\caption{Coarse and fine subdomain solution errors at different mesh }
\begin{center}
\begin{tabular}{ |c|c|c|c|c|c| } 
 \hline
 $h_{c}$ & $h_{f}$ & $||p-p_{h}||_{\Omega_{c}\times J}$ & $||p-p_{h}||_{\Omega_{f}\times J}$ & DOF & CPUTIM \\ 
 \hline
 1/10 & 1/40 & 0.0945 & 0.0159 & 11080 & 2.22\\ 
 1/20 & 1/80 & 0.0646 & 0.0096 & 88640 & 13.29\\
  1/40 & 1/160 & 0.0555 & 0.0084 & 709120 & 292.59 \\ 
 \hline
\end{tabular}
\end{center}
 \label{tab:err}
\end{table}
Table \ref{tab:err} shows the coarse and fine subdomain solution errors for increasing mesh refinement. Here, $h_{f}$ and $h_{c}$ are the fine and coarse subdomain h-refinements, respectively in both spatial and temporal dimensions. As expected, for this space-time monolithic, linear system at hand only one linear solve is required to reach the solution. Table \ref{tab:err} also shows the degrees of freedom associated with each of the three cases and the corresponding CPU time. Figure \ref{fig:manu} show the pressure solutions for the linear parabolic problem at the three refinement levels.
\section{Numerical results}\label{sec:results}
In this section we present two numerical experiments for practical problems of interest in the sub-surface porous medium (1) single phase flow and tracer transport, and (2) non-linear multiphase flow. These problems are chosen to demonstrate the applicability of the space-time domain decomposition approach, proposed in this work, to a wide range of subsurface physical processes. The single phase flow and tracer transport introduces a mild non-linearity in the form of slightly compressible fluid description whereas the two-phase flow problem is highly non-linear due to the description of relative permeability and capillary pressure functionals in the model formulation. Note that these two particular numerical experiments were chosen only to demonstrate and confirm expected resolution of non-linear convergence issues and does not preclude applicability to a more general class of problems.

\subsection{Single Phase Flow and Transport Problem}
We first present a numerical experiment for a slightly compressible, single phase flow and non-reactive tracer transport in porous medium. The transported component is assumed to be a tracer suchthat the fluid density remains invariant with changes in tracer concentration. The model formulation for a single phase, slightly compressible flow and tracer transport in porous medium is given by,
\begin{eqnarray}
\label{eqn:1phcon}
\frac{\partial}{\partial t} \left(\phi \rho\right) + \nabla \cdot \bs{u} = q \quad \text{in } \Omega\times J,\\
\label{eqn:1phtra}
\frac{\partial}{\partial t} \left(\phi \rho c\right) + \nabla \cdot \left( \bs{u} c +  \bs{z}\right) = q\hat{c}\quad \text{in } \Omega\times J,
\end{eqnarray}
with the Darcy velocity $u$ and diffusive flux $j$ defined as,
\begin{eqnarray}
\label{eqn:1phdar}
\bs{u} = -\frac{K}{\mu}\rho\nabla p \quad \text{in } \Omega\times J, and\\
\label{eqn:1phdif}
\bs{z} = -\phi \rho D\nabla c \quad \text{in } \Omega\times J
\end{eqnarray}
, respectively. Further the boundary and initial conditions are given by,
\begin{eqnarray}
\bs{u}\cdot \bs{\nu} = 0,\quad \bs{z}\cdot \bs{\nu} = 0 \quad \text{on } \partial\Omega\times J,\text{ and}
\end{eqnarray}
\begin{eqnarray}
p = p^{0}, \quad c = c^{0}  \quad \text{at } \partial\Omega\times \{0\}
\end{eqnarray}
, respectively. Here, $\bs{u}$ and $\bs{z}$ are the Darcy and diffusive fluxes, respectively, $\rho$ the density, $c$ and $p$ the pressure and concentration unknowns, respectively, $D$ is the constant (spatial and temporal) diffusion coefficient, $q$ the source/sink term, $\hat{c}$ the injection concentration at the source, $\nu$ the unit outward normal, $K$ is the spatially varying diagonal permeability tensor, $\mu$ the viscosity, $\phi$ the porosity, and $p^{0}$ and $c^{0}$ the initial conditions for pressure and concentration, respectively. The density is defined as a non-linear function of pressure given by, 
\begin{equation}
\rho = \rho_{ref}e^{c(p-p_{ref})}.
\label{eqn:sliden}
\end{equation}
Here, $\rho_{ref}$ is the reference density at the reference pressure $p_{ref}$, and $c$ is the fluid compressibility. The single phase flow and tracer transport model formulation presented here assumes a constant molecular diffusion for ease of description of a fully discrete mixed formulation. However, a spatially and temporally varying diagonal tensor can also be considered as the diffusion coefficient to account for hydrodynamic dispersion. A space-time, fully discrete form of the above model formulation is described in \ref{apx:a}. The fluid compressibility, viscosity, and density are taken to be 10$^{-6}~psi^{-1}$,1cP, and 64$lb/ft^{3}$,  respectively.  

The computational domain is (110ft $\times$ 30ft $\times$ 1ft $\times$ 100days) with the fine and coarse subdomains discretized using grid elements of size 0.5ft $\times$ 0.5ft $\times$ 1ft$\times$ 1day and 5ft $\times$ 5ft $\times$ 1ft $\times$ 5days, respectively. The fine subdomain is refined by a factor of 10 and 5 times with respect to the coarse subdomain in space and time, respectively. A source and sink term is considered at the bottom left and top right corners of the domain using an rate specified injection and pressure specified production wells, respectively. The injection well is water-rate specified at 4 STB/day with a non-reactive tracer injection at concentration 1. The production well is pressure specified at 1000 $psi$. A homogeneous porosity and permeability of 0.2 and 50 mD, respectively are assumed for the entire domain with a constant diffusion coefficient of 0.1. Further, the initial reservoir pressure and concentration are taken to be 1000 $psi$ and $0.0$, respectively. 

\begin{figure}[H]
\begin{center}
\includegraphics[width=7.5cm,trim=3cm 4.5cm 2cm 4cm, clip]{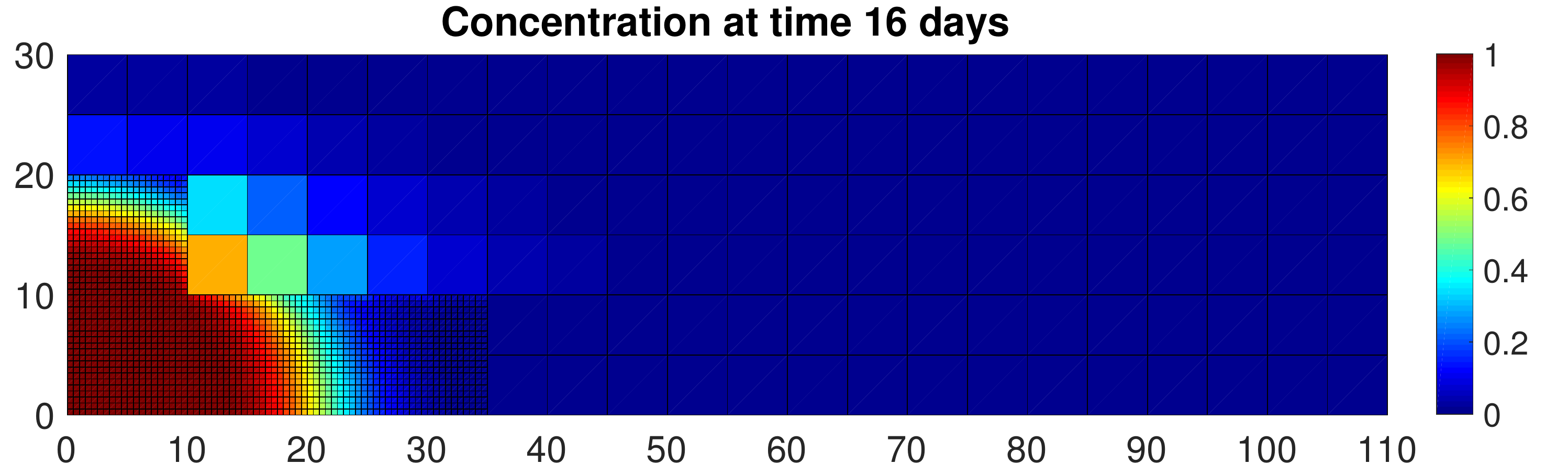}
\includegraphics[width=7.5cm,trim=3cm 4.5cm 2cm 4cm, clip]{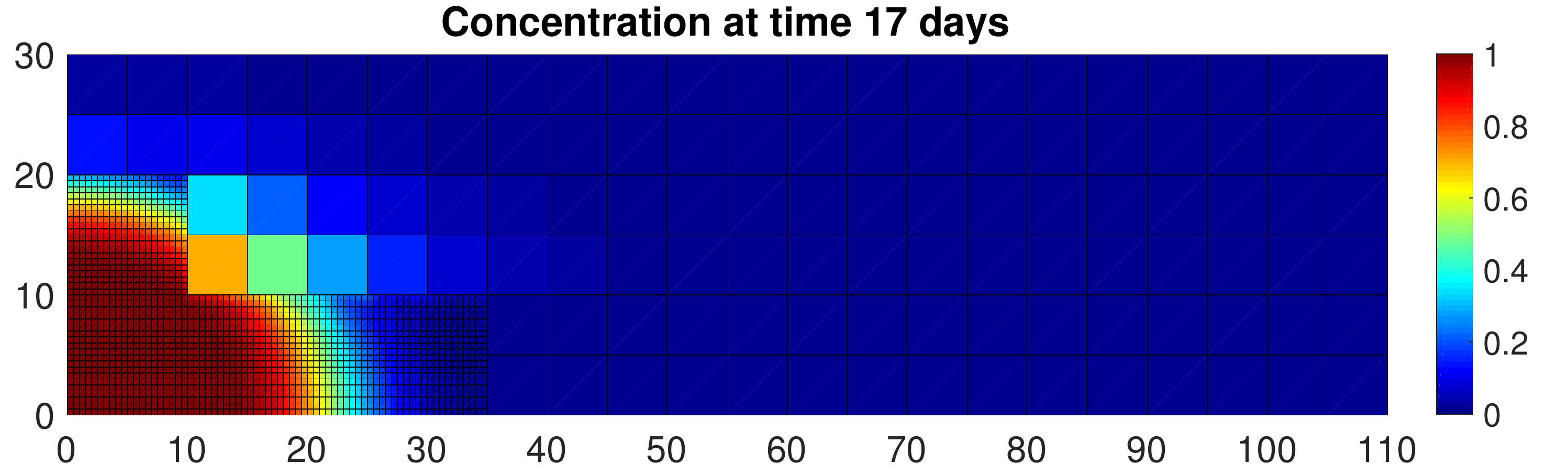}
\includegraphics[width=7.5cm,trim=3cm 4.5cm 2cm 4cm, clip]{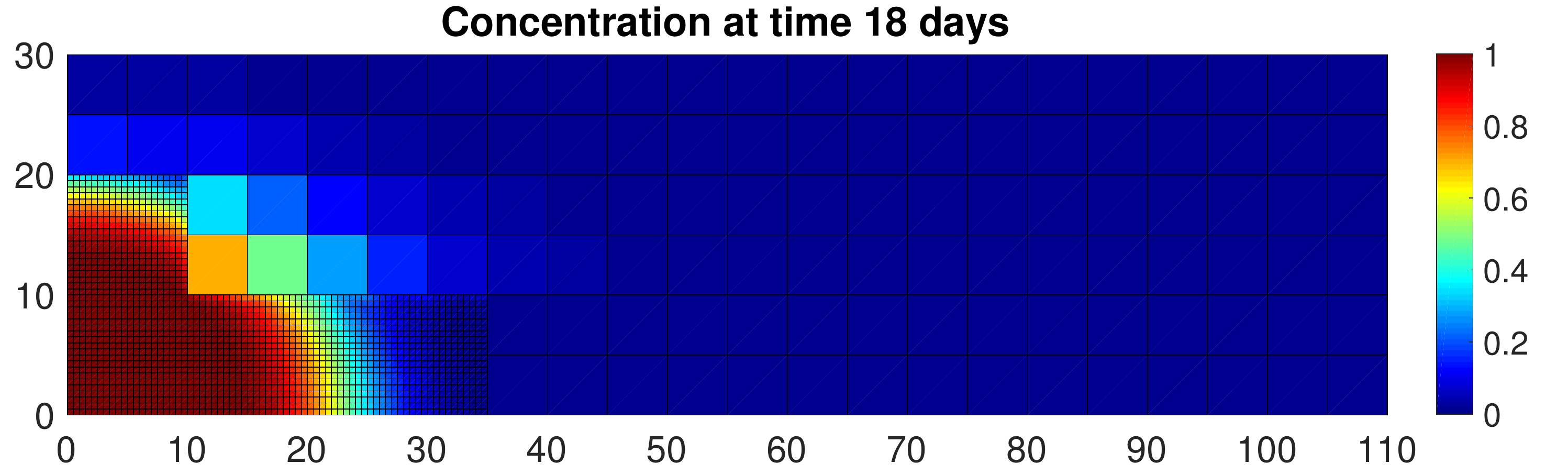}
\includegraphics[width=7.5cm,trim=3cm 4.5cm 2cm 4cm, clip]{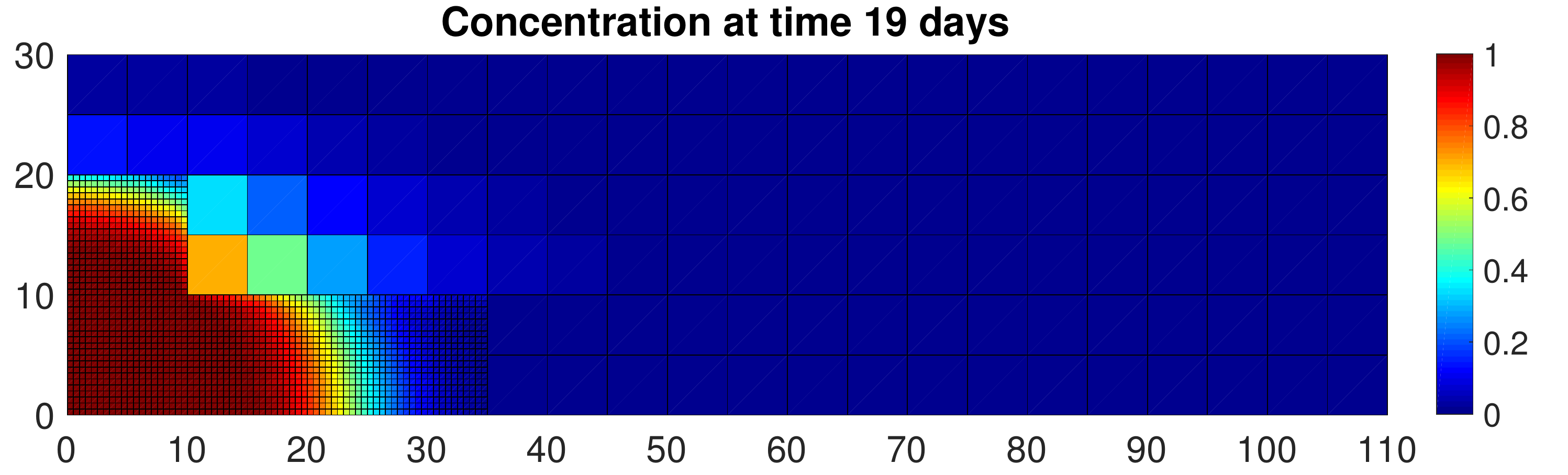}
\includegraphics[width=7.5cm,trim=3cm 4.5cm 2cm 4cm, clip]{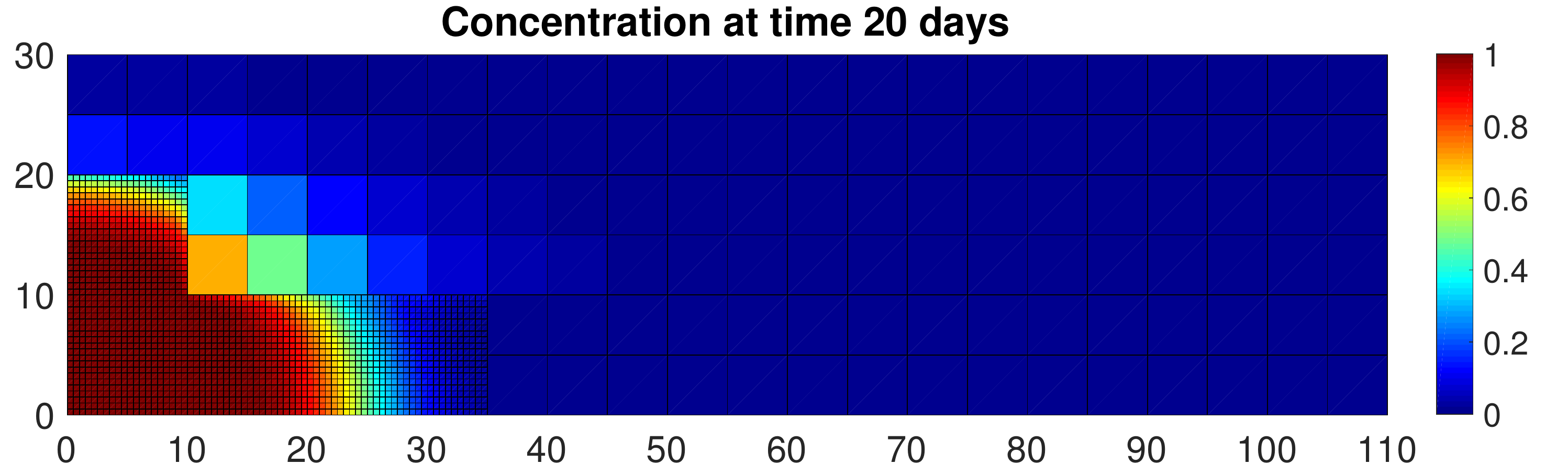}
\includegraphics[width=7.5cm,trim=3cm 4.5cm 2cm 4cm, clip]{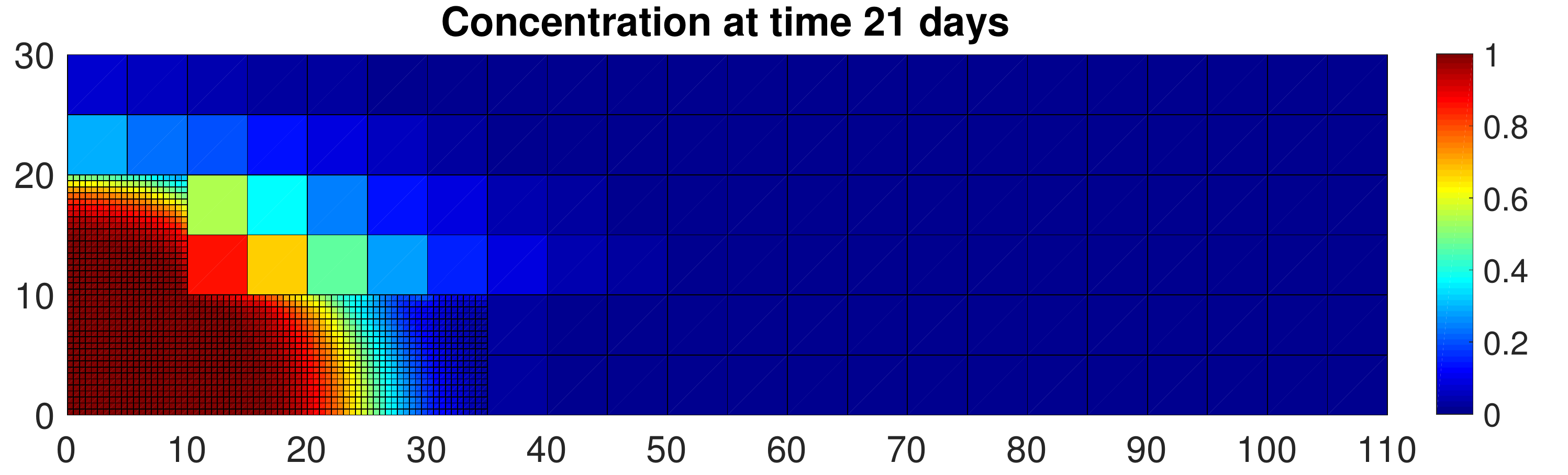}
\caption{Concentration distribution at different times}
\label{fig:1phconc}
\end{center}
\end{figure}

Figure \ref{fig:1phconc} shows the evolution of concentration distribution with time over the spatial domain. The concentration in the fine subdomain changes faster with small time-step increments of 1 day whereas the changes in the coarse subdomain occur every coarse time-step increment of 5 days. The fine and coarse subdomains are chosen in order to easily demonstrate these differences and this choice not restrictive for the space-time domain decomposition approach proposed here.

\subsection{Two Phase Flow Problem}
We describe the immiscible, two-phase, slightly compressible flow in porous medium model formulation with oil and water phase mass conservation and constitutive equations along with the boundary and initial conditions. The fully discrete space-time formulation, similar to single phase slightly compressible flow in the previous subsection, is presented in \ref{apx:b}. Further details regarding the two-phase flow model formulation can be found in \cite{singhapprox}.
The mass conservation equation for phase $\alpha$ is given by,
\begin{equation}
\frac{\partial\left(\phi\rho_{\alpha}s_{\alpha}\right)}{\partial t} + \div \u_{\alpha} = q_{\alpha} \text{~in~} \Omega \times J,
\label{eqn:2phcon}
\end{equation}
where $\phi$ and $K$ have their usual meanings as described before, and $\rho_{\alpha}$, $s_{\alpha}$, $u_{\alpha}$ and $q_{\alpha}$ are density, saturation, velocity and source/sink term, respectively of phase $\alpha$. The constitutive equation for the corresponding phase $\alpha$ is given by Darcy's law as,
\begin{equation}
\u_{\alpha} = -K\rho_{\alpha}\frac{k_{r\alpha}}{\mu_{\alpha}}\left(\grad p_\alpha-\rho_{\alpha}\g\right) \text{~in~} \Omega \times J
\end{equation}
Further, $k_{r\alpha}$, $\mu_{\alpha}$ and $p_{\alpha}$ are the relative permeability, viscosity and pressure of phase $\alpha$. 
Although not restrictive, for the sake of simplicity we assume no flow boundary conditions.
\begin{eqnarray}
\u_{\alpha} \cdot \bs{\nu} = 0 \text{~on~} \partial \Omega \times J\\
p_{\alpha} = p_{\alpha}^{0}, \quad s_{\alpha} = s_{\alpha}^{0},  \text{~at~} \Omega \times \{0\}
\end{eqnarray}
Here, $p^{0}_{\alpha}$, $s_{\alpha}^{0}$ are the initial conditions for pressure and saturation of phase $\alpha$.
Furthermore, the phase saturations $s_{\alpha}$ obey the constraint,
\begin{equation}
\sum_{\alpha}s_{\alpha} = 1.
\label{eqn:2phsat}
\end{equation}
We assume capillary pressure and relative permeabilities to be continuous and monotonic functions of phase saturations, 
\begin{equation}
p_{c}=f(s_{o}) = p_{w}-p_{o},
\label{eqn:2phcap}
\end{equation}
\begin{equation}
k_{r\alpha} = k_{r\alpha}(s_{\alpha}).
\end{equation}
The oil and water phase are assumed to slightly compressible with phase densities evaluated using,
\begin{equation}
\rho_{\alpha} = \rho_{\alpha,ref}\exp\left[c_{f\alpha}(p_{\alpha}-p_{\alpha,ref})\right].
\end{equation}
Here, $c_{f\alpha}$ is the compressibility and $\rho_{\alpha,ref}$ is the density of phase $\alpha$ at the reference pressure $p_{\alpha,ref}$. 

\begin{figure}[H]
\begin{center}
\includegraphics[width=7cm,trim=0cm 0cm 0cm 0cm, clip]{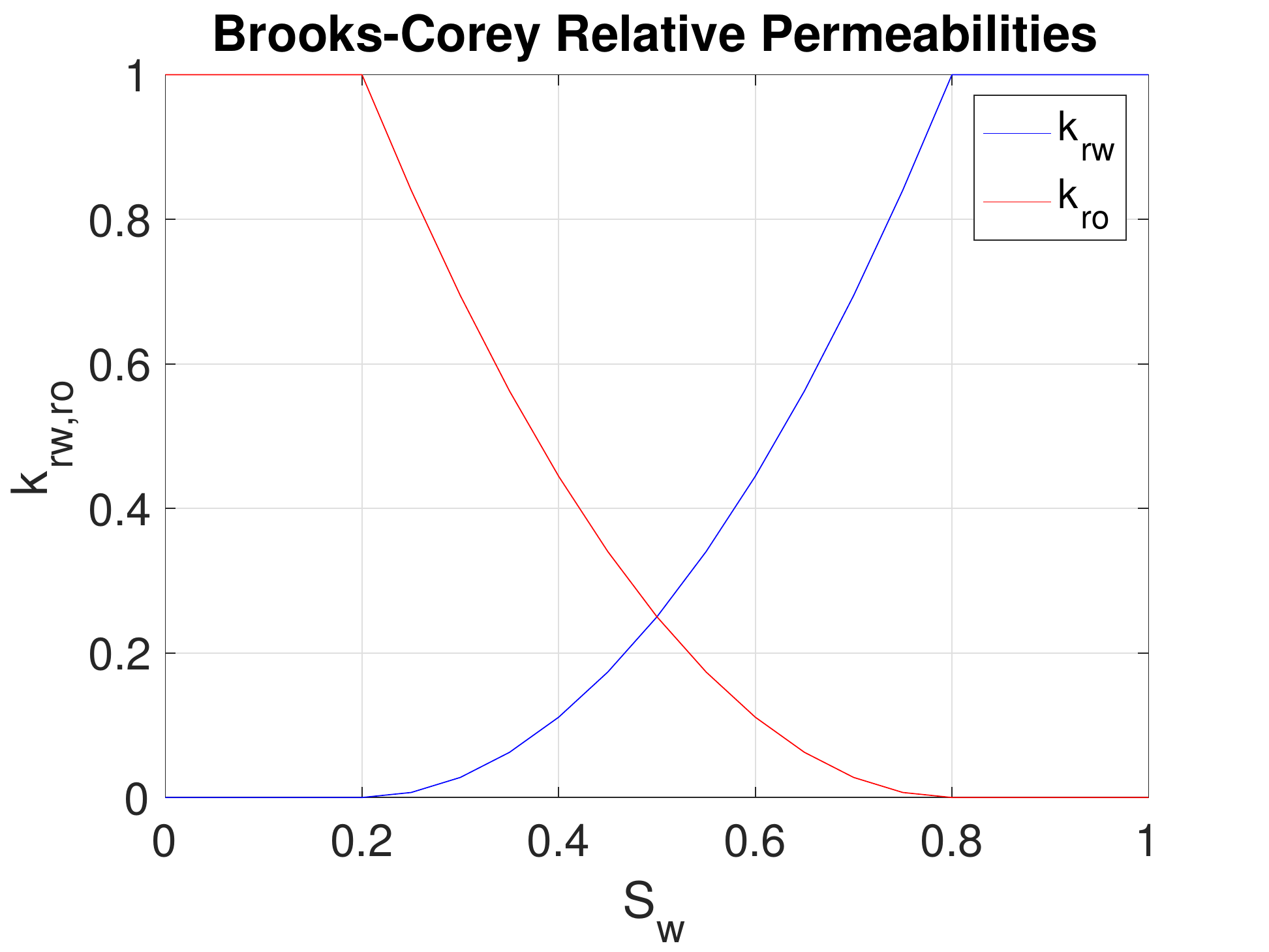}
\includegraphics[width=7cm,trim=0cm 0cm 0cm 0cm, clip]{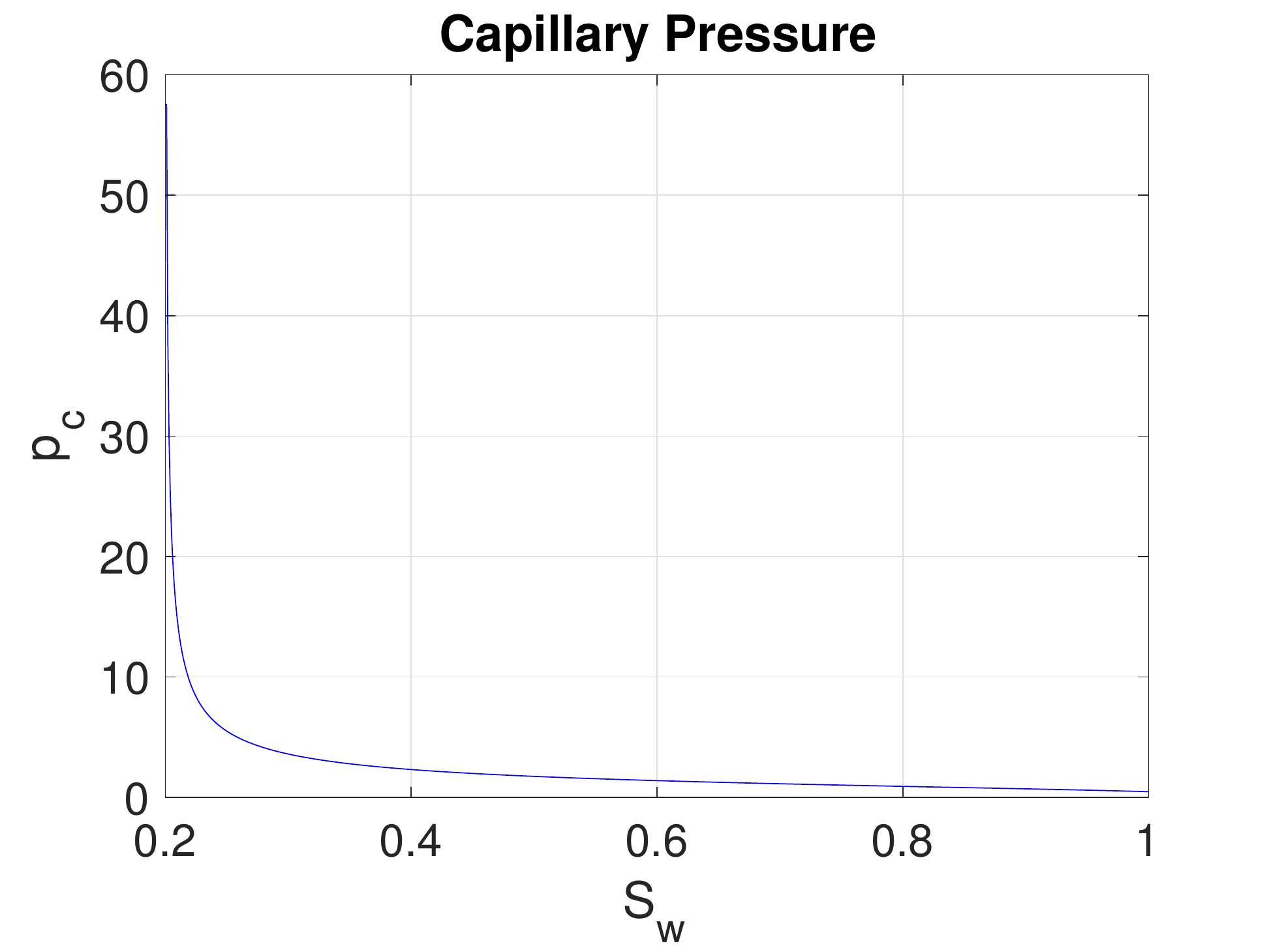}
\caption{Relative permeability (left) and capillary pressure (right) curves for the two-phase flow problem}
\label{fig:relcap}
\end{center}
\end{figure}

\begin{figure}[H]
\begin{center}
\includegraphics[width=7.5cm,trim=2.5cm 6cm 2cm 5cm, clip]{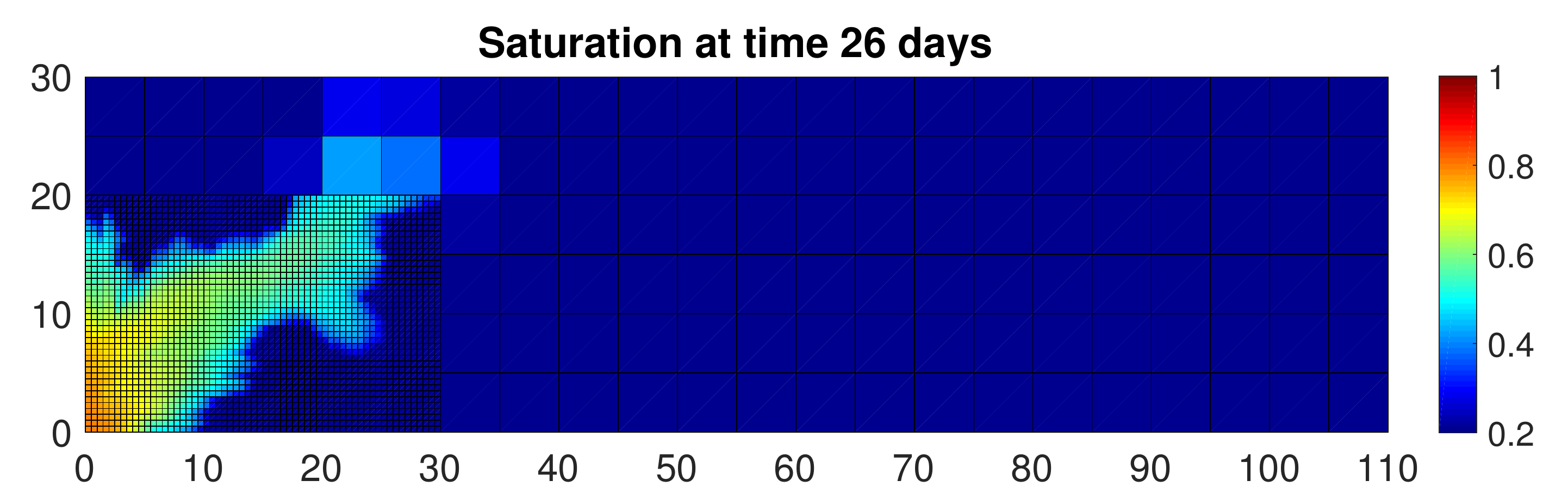}
\includegraphics[width=7.5cm,trim=2.5cm 6cm 2cm 5cm, clip]{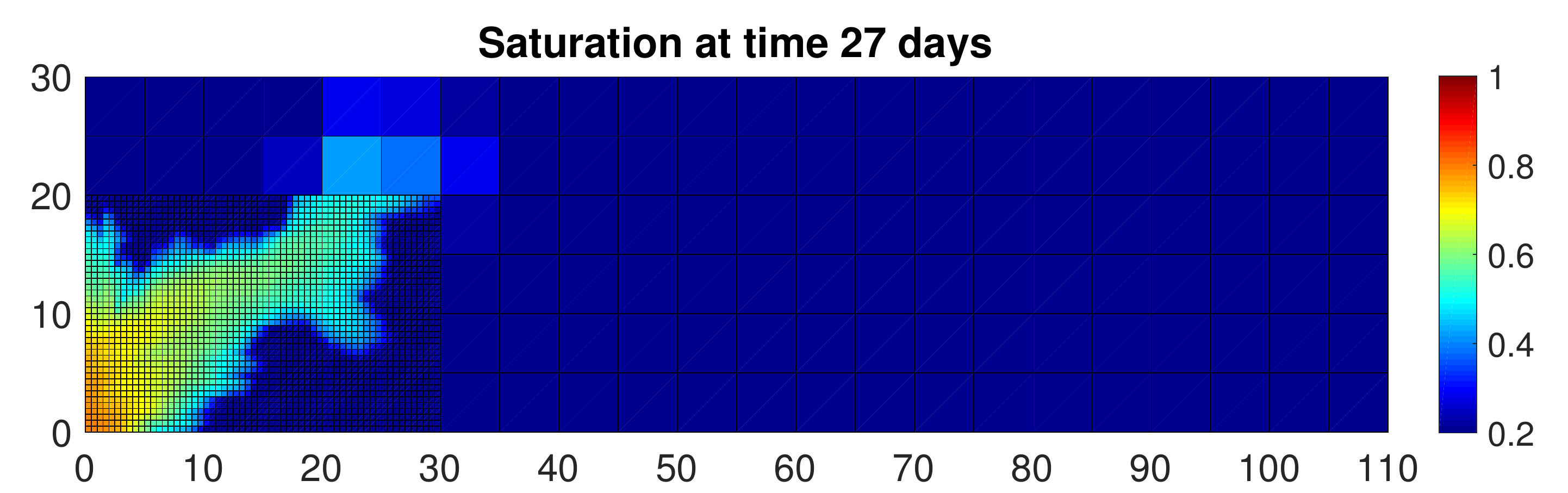}
\includegraphics[width=7.5cm,trim=2.5cm 6cm 2cm 5cm, clip]{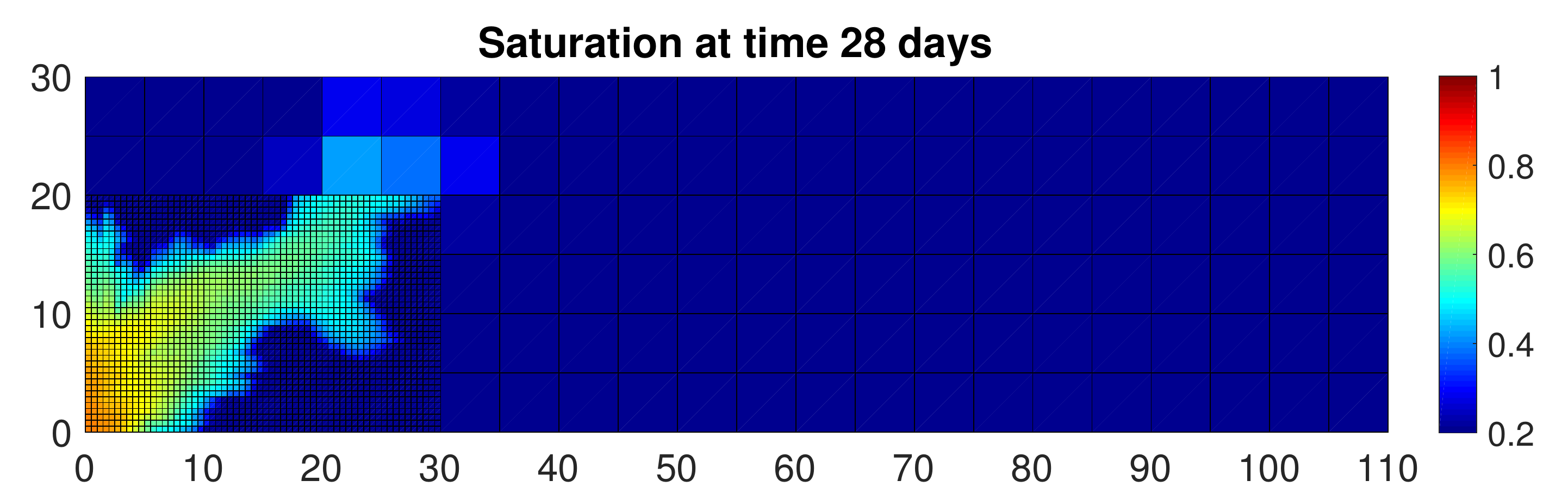}
\includegraphics[width=7.5cm,trim=2.5cm 6cm 2cm 5cm, clip]{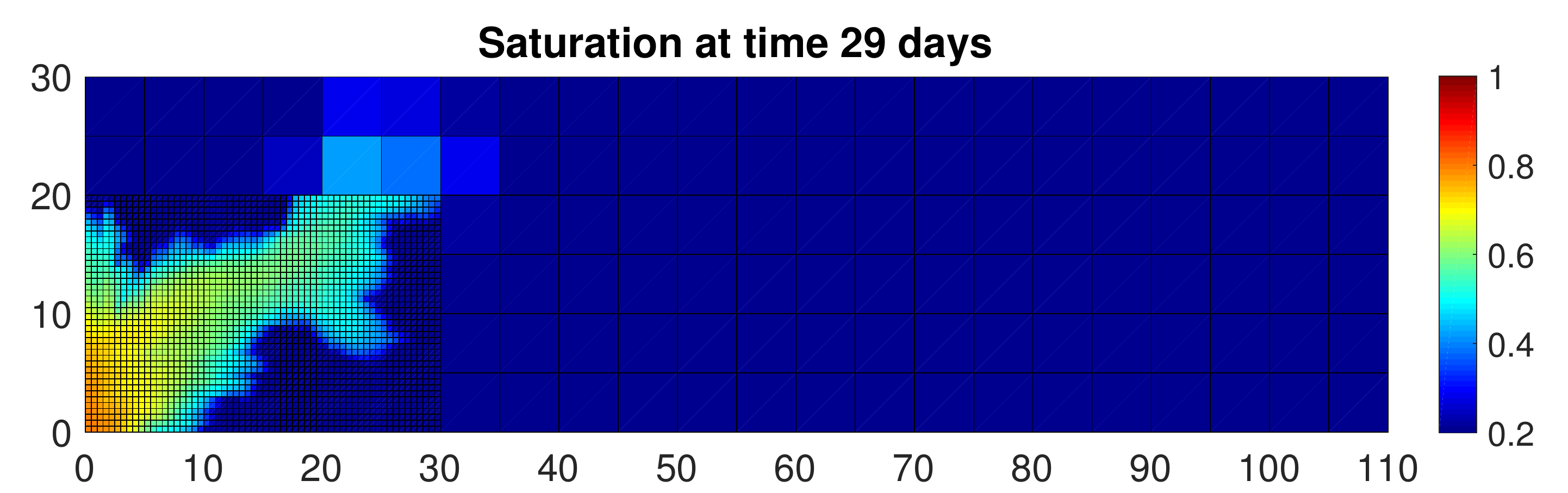}
\includegraphics[width=7.5cm,trim=2.5cm 6cm 2cm 5cm, clip]{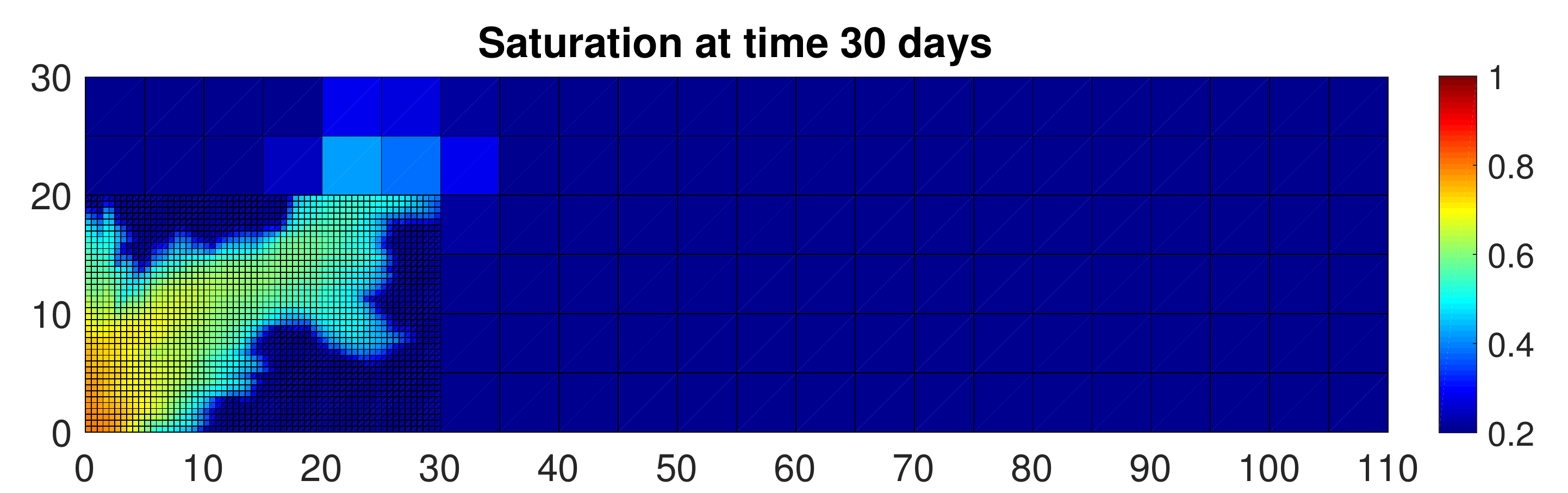}
\includegraphics[width=7.5cm,trim=2.5cm 6cm 2cm 5cm, clip]{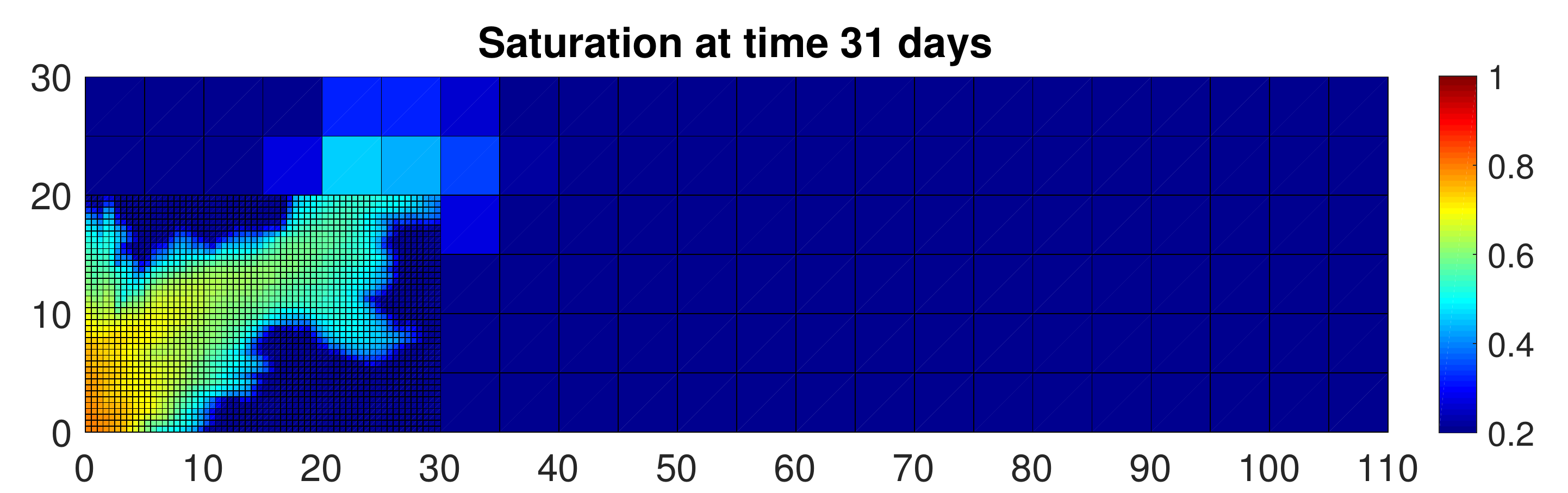}
\caption{Saturation distribution at different times}
\label{fig:2phsat}
\end{center}
\end{figure}

For the numerical experiment, the computational domain is kept the same as before (110ft $\times$ 30ft $\times$ 1ft $\times$ 40days) with the fine and coarse subdomains discretized using grid elements of size 0.5ft $\times$ 0.5ft $\times$ 1ft $\times$ 1day and 5ft $\times$ 5ft $\times$ 1ft $\times$ 5days, respectively. The fine subdomain is refined by a factor of 10 and 5 times with respect to the coarse subdomain in space and time, respectively. 
\begin{figure}[H]
\begin{center}
\includegraphics[width=7.5cm,trim=2.5cm 4.5cm 2.5cm 4cm, clip]{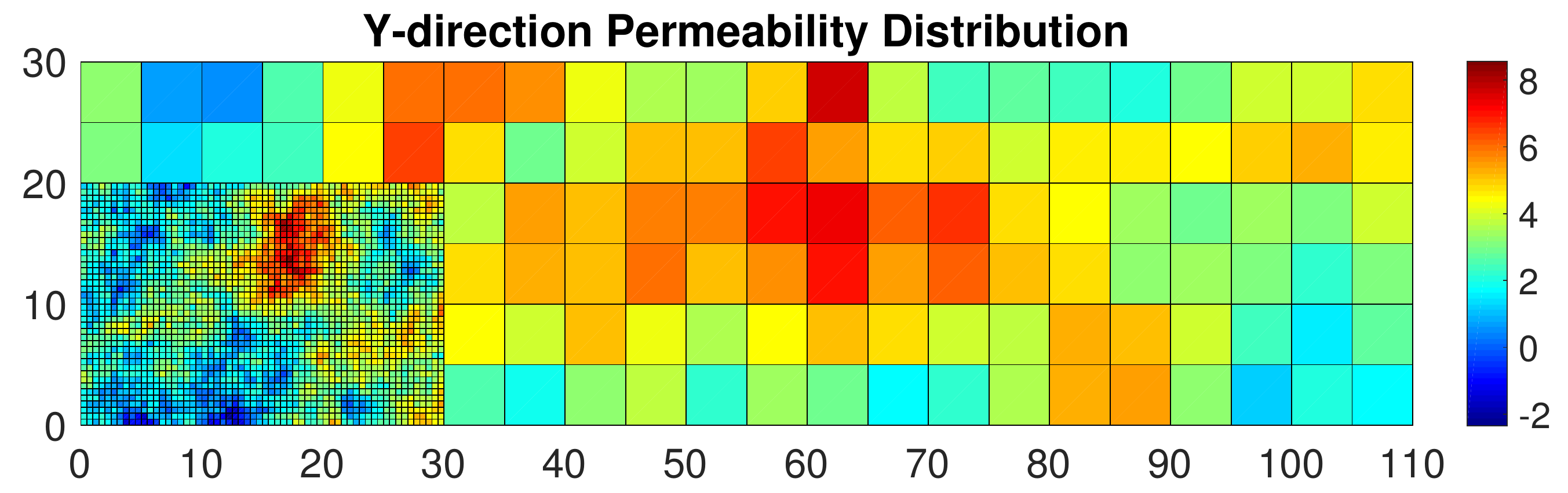}
\includegraphics[width=7.5cm,trim=2.5cm 4.5cm 2.5cm 4cm, clip]{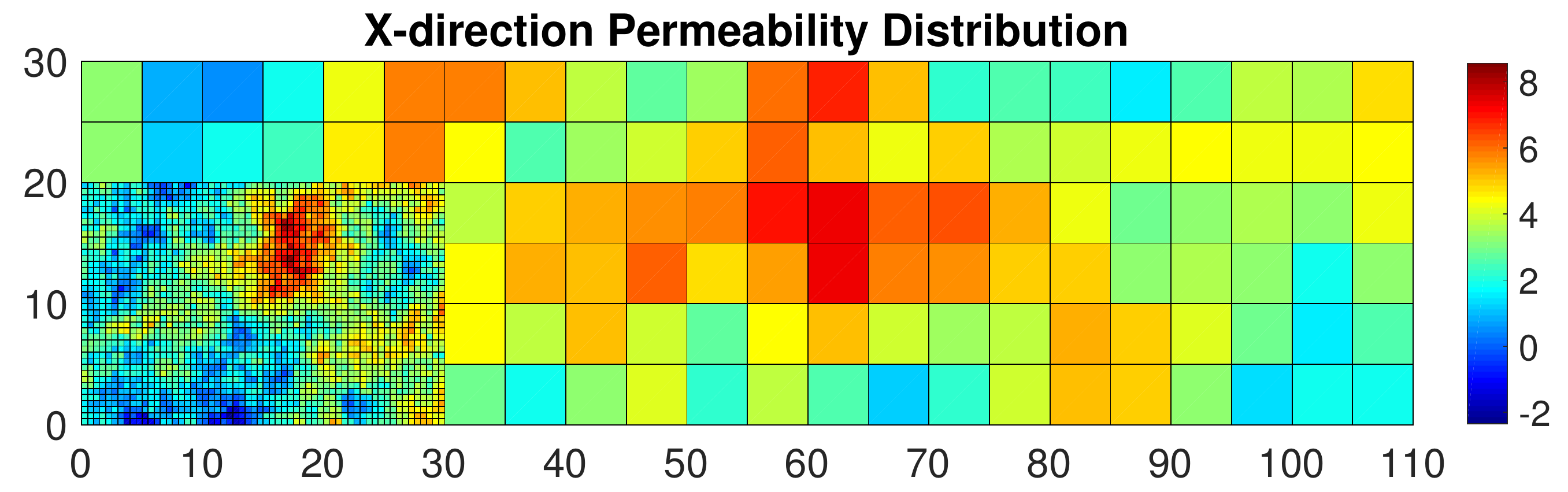}
\caption{Log-scale permeability distribution for the horizontal (left) and vertical direction (right)}
\label{fig:perm}
\end{center}
\end{figure}
The fluid and reservoir properties are adapted from the SPE10 \citep{spe10} dataset with an assumed homogeneous, spatial distribution for porosity of 0.2. The log (natural log) scale, permeability distribution in the horizontal (y) and vertical (x) directions are shown in Figure \ref{fig:perm}. The oil and water phase compressibilities are taken to be 1$\times$10$^{-4}$ and 3$\times$10$^{-6}~psi^{-1}$ , respectively, and densities 53$lb/ft^{3}$ and 64$lb/ft^{3}$, respectively. Further, the fluid viscosities are assumed to be 3 and 1 $cP$ for the oil and water phases, respectively. Additionally, a Brook's Corey model \cite{brooks}, Eqn. \eqref{eqn:brooks}, is considered for the two-phase relative permeabilities with endpoints s$_{or}$ = s$_{wirr}$ = 0.2 and $k^{0}_{ro}=k^{0}_{rw}=1.0$, and model exponents $n_{o}=n_{w}=2$. Figure \ref{fig:relcap} shows the relative permeability and capillary pressure curves as functions of saturation.
\begin{equation}
\begin{aligned}
k_{rw} &= k^{0}_{rw}\left(\frac{s_{w}-s_{wirr}}{(1-s_{or}-s_{wirr}}\right)^{n_{w}}\\
k_{ro} &= k^{0}_{ro}\left(\frac{s_{o}-s_{or}}{(1-s_{or}-s_{wirr}}\right)^{n_{o}}
\label{eqn:brooks}
\end{aligned}
\end{equation}
The capillary pressure function is defined using the van Genuchten model \cite{genuchten} given by Eqn. \eqref{eqn:genuchten}. The model parameters a, b, and c are chosen to be 0.8 psi, 0.6255, and 2.67, respectively.
\begin{equation}
p_{c}(s_{w}) = a \left[ (s_{w}-s_{wirr})^{-1/b} - 1\right]^{1/c}
\label{eqn:genuchten}
\end{equation}
As in the previous numerical experiment, a source and a sink are considered  at the bottom left and top right corners of the domain using an rate specified injection and pressure specified production wells, respectively. The injection well is water-rate specified at 1 STB/day whereas the production well is pressure specified at 1000 $psi$. Further, the initial reservoir pressure and saturation are taken to be 1000 $psi$ and $0.2$, respectively. Figure \ref{fig:2phsat} shows the evolution of the saturation front with time over the entire domain. The simulation was ran for a total of 40 days however we show the saturation distribution starting from 26 days up until 31 to demonstrate faster changes occurring in the fine domain compared to the coarse subdomain. The saturation distribution changes in the coarse subdomain at every coarse time-step increment (5 days) as opposed to fine subdomain where it changes at every fine time-step increment (1 day). 
The choice of the fine subdomain, at the bottom left part of the domain, is made to capture the evolution of the saturation front starting from the injection well. It is easy to see that the non-linear functions of saturation (such as relative permeability and capillary pressure) manifest highly non-linear behavior in the region where saturation changes are large. Here, non-linear solvers such as Newton-Raphson method are marred by small time-step sizes.  An increase in time-step size often results in an increase in non-linear iterations or convergence issues.  This is further exacerbated by mesh refinement (local or global) and consequently increased computational costs. A choice of fine subdomain (space and time) in the vicinity of the saturation front allows us to not only to gain accuracy but also circumvent convergence issues associated with the non-linear solver. The computational savings are self evident since small time-step increments are only necessary in the fine subdomain as opposed to the entire domain. In fact, for the linear flow and transport problem, in section \ref{sec:numcon}, the computational cost is directly related to the space-time degrees of freedom.

\section{Conclusions}
A space-time, domain decomposition approach is presented for flow and transport problems in subsurface porous medium using an enhanced velocity approach. We also described an efficient space-time monolithic solver for this approach that does not require subdomain iterations for computational efficiency. In fact, any linear problem requires only one linear solve to reach the solution, as expected. A preliminary numerical convergence analysis for a linear parabolic problem allows us to confirm that the proposed approach converges as the mesh is refined in both space and time. A rigorous derivation of the a-priori error estimates equipped with appropriate norms in order to determine exact convergence rates is reserved for a future work. We presented two, non-linear flow and transport problems differing in the degree of non-linearity to demonstrate the general applicability of this space-time, domain decomposition approach to a wide range of problems in subsurface porous medium. The two-phase, non-linear flow problem with space-time, mesh refinement in the vicinity of injection well was able to circumvent non-linear solver convergence issues without imposing restrictively small time-step sizes over the entire computational domain. Further, the proposed domain decomposition approach along with the numerical solution algorithm renders us a massively parallel, concurrent in time, computational framework allowing us to use parallel in time, linear solvers and preconditioners without compromising computational efficiency.

\appendix
\section{Appendix A: Fully Discrete Single Phase Flow and Transport Formulation} \label{apx:a}
The mixed variational form of Eqns. \eqref{eqn:1phcon} thru \eqref{eqn:1phdif} is: find $\bs{u}_{h}^{t} \in \bs{V}_{h}^{t,*}$,  $\bs{z}_{h}^{t} \in \bs{V}_{h}^{t,*}$, $c_{h}^{t} \in W_{h}^{t}$, and $p_{h}^{t} \in W_{h}^{t}$ such that,
\begin{eqnarray}
\label{eqn:wcon}
\left( \frac{\partial}{\partial t} \phi \rho, w \right)_{\Omega \times J} + \left(\nabla \cdot \bs{u}_{h}^{t}, w \right)_{\Omega \times J} = \left(q,w\right)_{\Omega \times J} \quad \bs{w}\in\bs{W}_{h}^{t},\\
\label{eqn:wtra}
\left( \frac{\partial}{\partial t} \phi \rho c_{h}^{t}, w \right)_{\Omega \times J} + \left(\nabla \cdot \left( \bs{u}_{h}^{t} c_{h}^{t}\right), w \right)  - \left(\nabla \cdot \bs{z}_{h}^{t}, w \right)_{\Omega \times J} = \left(q\hat{c},w\right)_{\Omega \times J} \quad \bs{w}\in\bs{W}_{h}^{t},
\end{eqnarray}
\begin{eqnarray}
\label{eqn:wdar}
\left(\frac{\mu}{\rho}K^{-1} \bs{u}_{h}^{t}, \bs{v}\right)_{\Omega \times J} - \left(p_{h}^{t}, \nabla \cdot \bs{v}\right)_{\Omega \times J} = 0 \quad \bs{v}\in\bs{V}_{h}^{t,*},\\
\label{eqn:wflu}
\left(\frac{1}{\phi\rho}D^{-1} \bs{z}_{h}^{t}, \bs{v}\right)_{\Omega \times J} - \left(c_{h}^{t}, \nabla \cdot \bs{v}\right)_{\Omega \times J} = 0 \quad \bs{v}\in\bs{V}_{h}^{t,*}
\end{eqnarray}
Here, $\rho = \rho (p_{h}^{t})$. Similar to the linear, single phase flow description in Section \ref{sec:discrete} the solution can be written as,
\begin{equation}
p = \sum_{m=1}^{q}\sum_{i = 1}^{r} P_{i}^{m} w_{i}^{m}, \quad \bs{u} = \sum_{m=1}^{q}\sum_{i = 1}^{r+1} U_{i+\frac{1}{2}}^{m} \varphi_{i+\frac{1}{2}}^{m},
\end{equation}

\begin{equation}
c = \sum_{m=1}^{q}\sum_{i = 1}^{r} C_{i}^{m} w_{i}^{m}, \text{ and} \quad \bs{z} = \sum_{m=1}^{q}\sum_{i = 1}^{r+1} Z_{i+\frac{1}{2}}^{m} \varphi_{i+\frac{1}{2}}^{m},
\end{equation}
using basis function as described earlier. We construct a non-linear, algebraic system of equations by testing the variational forms of the discrete constitutive (Darcy and diffusive flux) and conservation (phase and component mass) equations with $w_{j}^{n}$ and $\varphi_{j+\frac{1}{2}}^{n}$, respectively. We only describe integral terms in the weak variational form of the non-linear (slightly compressible), single phase flow model that differ from the description in Section \ref{sec:discrete}. The first term in the phase mass conservation Eqn. \eqref{eqn:wcon} can be expanded as,
\begin{equation}
\begin{aligned}
\left( \frac{\partial}{\partial t} \phi \rho\left(\sum_{m=1}^{q}\sum_{i=1}^{r} P_{i}^{m} w_{i}^{m}\right), w_{j}^{n} \right)_{\Omega \times J}   = &\left( \frac{\partial}{\partial t} \phi_{j}\rho\left(P_{j}^{n} w_{j}^{n}\right) , w_{j}^{n}\right) \\
& + \left( \phi_{j}\rho(P_{j}^{n})-\phi_{j}\rho(P_{j}^{n-1}), w_{j}^{n-1}\right)\\
= & \left( \phi_{j} \rho_{j}^{n} - \phi_{j} \rho_{j}^{n-1}\right) |E_{j}^{n-1}|.
\end{aligned}
\end{equation}
As before, for the coarse domain $(j_{0}+\frac{1}{2})^{+}$, the above integral remains unchanged. However, for the fine domain $(j_{0}+\frac{1}{2})^{-}$ we have,
\begin{equation}
\begin{aligned}
\left(\frac{\partial}{\partial t}\phi \rho , w_{j_{0}}^{n-\frac{2}{3}}\right)_{(\Omega \times J)} = & \left( \phi_{j_{0}}\rho_{j_{0}}^{n-\frac{2}{3}} - \phi_{j_{0}}\rho_{j_{0}}^{n-1}\right) |E_{j_{0}}^{n-1}|, \\
\left(\frac{\partial}{\partial t}\phi \rho , w_{j_{0}}^{n-\frac{1}{3}}\right)_{(\Omega \times J)} = & \left( \phi_{j_{0}}\rho_{j_{0}}^{n-\frac{1}{3}} -\phi_{j_{0}}\rho_{j_{0}}^{n-\frac{2}{3}}\right) |E_{j_{0}}^{n-\frac{2}{3}}|, \\
\left(\frac{\partial}{\partial t}\phi \rho , w_{j_{0}}^{n}\right)_{(\Omega \times J)} = & \left( \phi_{j_{0}}\rho_{j_{0}}^{n} - \phi_{j_{0}}\rho_{j_{0}}^{n-\frac{1}{3}}\right) |E_{j_{0}}^{n-\frac{1}{3}}|.
\end{aligned}
\end{equation}
Similarly, the accumulation (first) term in the transport equation Eqn. \eqref{eqn:wtra} is expanded for the coarse and fine subdomains as,
\begin{equation}
\left( \frac{\partial}{\partial t} \phi \rho c, w_{j_{0}+1}^{n} \right)_{\Omega \times J} =\left( \phi_{j_{0}+1} (\rho c)_{j_{0}+1}^{n}- \phi_{j_{0}+1} (\rho c)_{j_{0}+1}^{n-1}\right) |E_{j_{0}+1}^{n-1}|
\end{equation}
, and
\begin{equation}
\begin{aligned}
\left(\frac{\partial}{\partial t}\phi \rho c , w_{j_{0}}^{n-\frac{2}{3}}\right)_{(\Omega \times J)} = & \left( \phi_{j_{0}}\rho_{j_{0}}^{n-\frac{2}{3}} c_{j_{0}}^{n-\frac{2}{3}}- \phi_{j_{0}}\rho_{j_{0}}^{n-1} c_{j_{0}}^{n-1}\right) |E_{j_{0}}^{n-1}|, \\
\left(\frac{\partial}{\partial t}\phi \rho c , w_{j_{0}}^{n-\frac{1}{3}}\right)_{(\Omega \times J)} = & \left( \phi_{j_{0}}\rho_{j_{0}}^{n-\frac{1}{3}}c_{j_{0}}^{n-\frac{1}{3}} -\phi_{j_{0}}\rho_{j_{0}}^{n-\frac{2}{3}}c_{j_{0}}^{n-\frac{2}{3}}\right) |E_{j_{0}}^{n-\frac{2}{3}}|, \\
\left(\frac{\partial}{\partial t}\phi \rho c, w_{j_{0}}^{n}\right)_{(\Omega \times J)} = & \left( \phi_{j_{0}}\rho_{j_{0}}^{n}c_{j_{0}}^{n} - \phi_{j_{0}}\rho_{j_{0}}^{n-\frac{1}{3}}c_{j_{0}}^{n-\frac{1}{3}}\right) |E_{j_{0}}^{n-\frac{1}{3}}|
\end{aligned}
\end{equation}
, respectively. The advection or second term in Eqn. \eqref{eqn:wtra} is approximated for the fine and coarse subdomains with,
\begin{equation}
\begin{aligned}
\left(\nabla \cdot \bs{u}c, w_{j_{0}}^{n-\frac{2}{3}}\right)_{\Omega \times J} & \approx (U\tilde{c})_{j_{0}+\frac{1}{2}}^{n-\frac{2}{3}} - (U\tilde{c})_{j_{0}-\frac{1}{2}}^{n-\frac{2}{3}},\\
\left(\nabla \cdot \bs{u}c, w_{j_{0}}^{n-\frac{1}{3}}\right)_{\Omega \times J} & \approx (U\tilde{c})_{j_{0}+\frac{1}{2}}^{n-\frac{1}{3}} - (U\tilde{c})_{j_{0}-\frac{1}{2}}^{n-\frac{1}{3}},\\
\left(\nabla \cdot \bs{u}c, w_{j_{0}}^{n}\right)_{\Omega \times J} & \approx (U\tilde{c})_{j_{0}+\frac{1}{2}}^{n} - (U\tilde{c})_{j_{0}-\frac{1}{2}}^{n}
\end{aligned}
\label{eqn:lincon2f}
\end{equation}
, and
\begin{equation}
\begin{aligned}
\left(\nabla \cdot \bs{u}c, w_{j_{0}+1}^{n}\right)_{\Omega \times J} & \approx (U\tilde{c})_{j_{0}+\frac{3}{2}}^{n} - (U\tilde{c})_{j_{0}+\frac{1}{2}}^{n-\frac{2}{3}} - (U\tilde{c})_{j_{0}+\frac{1}{2}}^{n-\frac{1}{3}} - (U\tilde{c})_{j_{0}+\frac{1}{2}}^{n}\\
\end{aligned}
\label{eqn:lincon2c}
\end{equation}
, respectively. Here, $\tilde{c}$ is the upwinded concentration to add sufficient diffusion for stability of the numerical scheme and is defined as,
\begin{equation}
\tilde{c}^{n}_{j+\frac{1}{2}} =
\begin{cases}
{c}^{n}_{j}, \quad & \text{ if } U^{n}_{j+\frac{1}{2}}>0\text{ and }\\[3ex]
{c}^{n}_{j+1}, \quad & \text{ otherwise}.
\end{cases}
\end{equation}

It is easy to see that the second and third term in Eqns. \eqref{eqn:wcon} and \eqref{eqn:wtra}, respectively are similar can be evaluated using Eqns. \eqref{eqn:lincon2f} and \eqref{eqn:lincon2c} for the fine and coarse subdomains, respectively. Further, the right hand side in Eqns. \eqref{eqn:wcon} and \eqref{eqn:wtra} are also similar and can be evaluated as in Eqn. \eqref{eqn:linconsource}. For the constitutive equations (Darcy and diffusive flux), the first terms in Eqns. \eqref{eqn:wdar} and \eqref{eqn:wflu} are then expressed as,
\begin{equation}
\begin{aligned}
\left( \frac{\mu}{\rho}K^{-1}\bs{u},\varphi_{j+\frac{1}{2}}^{n}\right)_{\Omega\times J} & \approx  \left( \frac{\mu}{\rho}K^{-1}\sum_{m=1}^{q} \sum_{i=1}^{r+1} U_{i+\frac{1}{2}}^{m}\varphi_{i+\frac{1}{2}}^{m}, \varphi_{j+\frac{1}{2}}^{n} \right)_{TM}\\
& = \frac{\mu}{2 \left| e_{j+\frac{1}{2}}^{n} \right|}\, \frac{2}{\left(\rho_{j}^{n} + \rho_{j+1}^{n}\right)}\left( \frac{h_{j}}{K_{j}} + \frac{h_{j+1}}{K_{j+1}} \right) U_{j+\frac{1}{2}}^{n},
\end{aligned}
\end{equation}
and
\begin{equation}
\begin{aligned}
\left( \frac{1}{\phi\rho}D^{-1}\bs{z},\varphi_{j+\frac{1}{2}}^{n}\right)_{\Omega\times J} & \approx  \left( \frac{1}{\phi \rho}D^{-1}\sum_{m=1}^{q} \sum_{i=1}^{r+1} Z_{i+\frac{1}{2}}^{m}\varphi_{i+\frac{1}{2}}^{m}, \varphi_{j+\frac{1}{2}}^{n} \right)_{TM}\\
& = \frac{D^{-1}}{2 \left| e_{j+\frac{1}{2}}^{n} \right|}\, \frac{2}{\left(\rho_{j}^{n} + \rho_{j+1}^{n}\right)} \left( \frac{h_{j}}{\phi_{j}} + \frac{h_{j+1}}{\phi_{j+1}} \right) Z_{j+\frac{1}{2}}^{n},
\end{aligned}
\end{equation}
respectively. The second terms in Eqns. \eqref{eqn:wdar} and \eqref{eqn:wflu} can be written for the fine and coarse subdomains similar to Eqns. \eqref{eqn:linflu2f} and \eqref{eqn:linflu2c}, respectively. The above terms gathered together give us a non-linear, algebraic system of equations in pressure ($p$), concentration ($c$), and Darcy ($\bs{u}$) and diffusive ($\bs{z}$) flux unknowns. The flux unknowns are eliminated, after Newton linearization of the aforementioned non-linear algebraic equations, using first and second Schur complement of the linearized equations resulting in a reduced system in pressure and concentration unknowns only.  

\section{Appendix B: Fully Discrete Two Phase Flow Formulation} \label{apx:b}
The expanded mixed variational form of Eqns. \eqref{eqn:2phcon} thru \eqref{eqn:2phsat} is: find $\bs{u}_{\alpha,h}^{t} \in \bs{V}_{h}^{t,*}$,  $\tilde{\bs{u}}_{\alpha,h}^{t} \in \bs{V}_{h}^{t,*}$, $S_{w,h}^{t} \in W_{h}^{t}$, and $p_{o,h}^{t} \in W_{h}^{t}$ such that,

\begin{equation}
\left(\frac{\partial}{\partial t} \phi\left(\rho_{w}s^{t}_{w,h} + \rho_{o}(1-s^{t}_{w,h})\right),w\right) + \left(\div \left(\bs{u}^{t}_{w,h}+\bs{u}^{t}_{o,h}\right),w\right) = \left(q_{w}+q_{o},w\right)
\label{eqn:totcon}
\end{equation}

\begin{equation}
\left(\frac{\partial}{\partial t} \left(\phi \rho_{w}s^{t}_{w,h} \right),w\right) + \left(\div \bs{u}^{t}_{w,h},w\right) = \left(q_{w},w\right),
\label{eqn:watcon}
\end{equation}

\begin{equation}
\left(K^{-1}\tilde{\u}^{t}_{o,h},\bs{v}\right) - \left(p^{t}_{o,h},\div \bs{v}\right) =0,
\label{eqn:oildar}
\end{equation}
\begin{equation}
\begin{aligned}
\left(K^{-1}\tilde{\u}^{t}_{w,h},\bs{v}\right) - \left(p^{t}_{o,h},\div \bs{v}\right) = -\left(p_{c},\div \bs{v}\right) ,
\end{aligned}
\label{eqn:watdar}
\end{equation}
\begin{equation}
\begin{aligned}
\left( \u^{t}_{\alpha,h},\bs{v}\right) = \left(\lambda_{\alpha} \tilde{\u}^{t}_{\alpha,h},\bs{v}\right),
\end{aligned}
\label{eqn:expand}
\end{equation}
with $w\in W$ and $\bs{v}\in \bs{V}$. Please note that $s_{o}$ and $p_{w}$ are eliminated in the above formulation using the algebraic constraints Eqns. \eqref{eqn:2phsat} and \eqref{eqn:2phcap}, respectively. Further, $\lambda_{\alpha}$ is defined as the mobility of phase $\alpha$ as,
\begin{equation}
\lambda_{\alpha}= \frac{k_{r\alpha}\rho_{\alpha}}{\mu_{\alpha}},
\end{equation}
An expanded mixed formulation \cite{weiser85, malgo02}, with additional auxiliary phase fluxes $\tilde{\bs{u}}_{\alpha}$, is used to avoid inverting zero phase relative permeabilities ($k_{r\alpha}$). As before,  the solution can be written as,
\begin{equation}
p_{o} = \sum_{m=1}^{q}\sum_{i = 1}^{r} P_{i}^{m} w_{i}^{m}, \quad \bs{u}_{\alpha} = \sum_{m=1}^{q}\sum_{i = 1}^{r+1} U_{\alpha,i+\frac{1}{2}}^{m} \varphi_{i+\frac{1}{2}}^{m},
\end{equation}
\begin{equation}
s_{w} = \sum_{m=1}^{q}\sum_{i = 1}^{r} S_{w,i}^{m} w_{i}^{m}, \text{ and} \quad \bs{\tilde{u}}_{\alpha} = \sum_{m=1}^{q}\sum_{i = 1}^{r+1} \tilde{U}_{\alpha,i+\frac{1}{2}}^{m} \varphi_{i+\frac{1}{2}}^{m},
\end{equation}
It is easy to see that all the terms in Eqns. \eqref{eqn:totcon} thru \eqref{eqn:watdar} can be expanded as before in \ref{apx:a}. The first and second terms in Eqn. \eqref{eqn:expand} are approximated as,
\begin{equation}
\begin{aligned}
(\bs{u}_{\alpha}, \bs{v})  \approx (\bs{u}_{\alpha}, \bs{v})_{Q}
	 & = (\bs{u}_{\alpha}, \phi^{n}_{j+\frac{1}{2}})_{Q} \\[1ex] 
	& = \sum_{m=1}^{q}\sum_{i = 1}^{r+1} U^{m}_{\alpha,i+1/2} (\phi^{m}_{i+\frac{1}{2}}, \phi^{n}_{j+\frac{1}{2}})_{Q}
	   = \frac{h_j{} + h_{j+1}}{2\, \rvert e^{n}_{j+\frac{1}{2}} \rvert} \,U^{n}_{\alpha,j+\frac{1}{2}}, \text{and}
	    \label{eqn:aux_left_bf}		
\end{aligned}
\end{equation}
\begin{equation} \label{eqn:aux_right_bf}
(\lambda_{\alpha}\tilde{\bs{u}}_{\alpha}, \bs{v}) \approx (\lambda^{*}_{\alpha}\tilde{\bs{u}}_{\alpha},\bs{v})_{Q} 
 = \frac{h_{j}+ h_{j+1}}{2\, \rvert e^{n}_{j+\frac{1}{2}} \rvert} \, \lambda^{*,n}_{\alpha,j+\frac{1}{2}} \tilde{U}^{n}_{\alpha,j+\frac{1}{2}}
\end{equation}
, respectively. Here, is $Q$ the quadrature rule in Eqn. \eqref{eqn:quad}, $\lambda^{*,n}_{\alpha,j+\frac{1}{2}}$ is the upwind mobility defined as,
\begin{equation}\label{eqn:upw_mob}
\lambda^{*,n}_{\alpha,j+\frac{1}{2}} = \rho^{*,n}_{\alpha,j+\frac{1}{2}} \dfrac{k^{r\alpha,*}_{j+\frac{1}{2}}}{\mu_{\alpha}} =
\begin{cases}
\dfrac{1}{2\mu_{\alpha}} \left(\rho^{n}_{\alpha,j} + \rho^{n}_{\alpha,j+1}\right) k_{r\alpha}(S^{n}_{\alpha,j}), \quad & \text{~if ~} \tilde{U}^{n}_{\alpha,j+\frac{1}{2}} > 0,\\[3ex]
\dfrac{1}{2\mu_{\alpha}}  \left(\rho^{n}_{\alpha,j} + \rho^{n}_{\alpha,j+1}\right) k_{r\alpha}(S^{n}_{\alpha,j+1}),  & \text{~otherwise}.
\end{cases}
\end{equation}
This gives us a non-linear, algebraic system of equations in pressure ($p_{o}$), saturation ($s_{w}$), and Darcy ($\bs{u}_{\alpha}$) and auxiliary ($\tilde{\bs{u}}_{\alpha}$) flux unknowns. As in \ref{apx:a}, the flux unknowns are eliminated using multiple Schur complements of the linearized system of equations resulting in a reduced system in pressure and saturations unknowns only.

\section*{Acknowledgements}
The authors would like to thank Center for Subsurface Modeling industrial affiliates for their continued support. We would also like to thank Prof. Ivan Yotov (Department of Mathematics, University of Pittsburgh) for the discussions on the enhanced velocity mixed finite element method.

%
\bibliographystyle{plain} 
\bibliography{ref}

\end{document}